\documentclass[11pt,twoside,leqno]{aomamlt2e} 
  \pageno{1033}
\received{June 14, 2004}
\revised{February 16, 2005}
   \renewcommand{\thesection}{\arabic{section}}
\makeatletter
\renewcommand{\@seccntformat}[1]{\csname
the#1\endcsname.\hspace{0.5em}\setcounter{Subsec}{0}\setcounter{Subsubsec}{0}}\makeatother

 \newtheorem{thm}{Theorem} 
  
 \newtheorem{lem}{Lemma}

\DeclareMathOperator{\PLA}{PLA}
\DeclareMathOperator{\real}{Re}
\renewcommand{\theenumi}{\rm (\roman{enumi})}

\begin{document}
\currannalsline{164}{2006} 

\title{Analytic representation of functions and\\ a new quasi-analyticity
threshold}

 \acknowledgements{Research supported in part by the Israel Science Foundation.}
 
\twoauthors{Gady Kozma}{Alexander Olevski\u\i}

 \institution{Weizmann Institute of Science, Rehovot 76100, Israel\\
\email{gady.kozma@weizmann.ac.il}
\\
\vglue-9pt
Tel Aviv University, Tel Aviv, Israel
\\
\email{olevskii@post.tau.ac.il}}

 \shorttitle{Analytic representation}   

\centerline{\bf Abstract}
\vglue4pt
We characterize precisely the possible rate of decay of the anti-analytic
half of a trigonometric series converging to zero almost everywhere.
 
\vglue-20pt
\phantom{up}
\section{Introduction}
\vglue-8pt

1.1.\quad In 1916, D.\ E.\ Menshov constructed an example of a nontrivial trigonometric
series on the circle $\mathbb{T}$ \begin{equation}
\sum_{n=-\infty}^{\infty}c(n)e^{int}\label{eq:Null}\end{equation}
which converges to zero almost everywhere (a.e.). Such series are
called null-series. This result was the origin of the modern theory
of uniqueness in Fourier analysis, 
see \cite{Z59}, \cite{B64}, \cite{KL87}, \cite{KS94}.

Clearly for such a series $\sum|c(n)|^{2}=\infty$. A less trivial
observation is that a null series cannot be analytic, that is, involve
positive frequencies only. Indeed,   it would then follow by Abel's
theorem that the corresponding analytic function
 \begin{equation}
F(z)=\sum_{n\geq0}c(n)z^{n}\label{eq:zPLA}\end{equation}
has nontangential boundary values equal to zero a.e.~on the circle
$|z|=1$. Privalov's uniqueness theorem (see below in \S\ref{sub:uniqhist})
now shows that $F$ is identically zero. 

\vskip6pt {\it Definition}.
We say that a function $f$ on the circle $\mathbb{T}$ belongs to
$\PLA$ (which stands for Pointwise Limit of Analytic series) if it
admits a representation\begin{equation}
f(t)=\sum_{n\geq0}c(n)e^{int}\label{eq:defPLA}\end{equation}
by an a.e.~converging series.

The discussion above shows that such a representation is unique. Further,
for example, $e^{-int}$ is not in $\PLA$ for any $n>0$ since multiplying
with $z^{n}$ would lead to a contradiction to Privalov's theorem.

If $f$ is an $L^{2}$ function with positive Fourier spectrum, or
in other words, if it belongs to the Hardy space $H^{2}$, then it
is in $\PLA$ according to the Carleson convergence theorem. On the
other hand, we proved in \cite{KO03} that $L^{2}$ contains in addition
$\PLA$ functions which are not in $H^{2}$. The representation (\ref{eq:defPLA})
for such functions is ``nonclassical'' in the sense that it is
different from the Fourier expansion.

One should contrast this phenomenon against some results in the Riemannian
theory (see \cite[Chap.~11]{Z59}) which say that whenever a representation
by harmonics is unique then it is the Fourier one. Compare for examples
the Cantor theorem to the du Bois-Reymond theorem. In an explicit
form this principle was stated in \cite{P23}: If a function $f\in L^{1}(\mathbb{T})$
has a \emph{unique} pointwise decomposition (\ref{eq:Null}) outside
of some compact $K$ then it is the Fourier expansion of $f$. Again,
for analytic expansions (\ref{eq:defPLA}) this is not true.

\demo{{\rm 1.2}} Taking a function $f$ from the ``nonclassic'' part of $\PLA\cap\, L^{2}$
and subtracting from the representation (\ref{eq:defPLA}) the Fourier
expansion of $f$, one gets a null-series with a small anti-analytic
part in the sense that\[
\sum_{n<0}|c(n)|^{2}<\infty.\]
Note that there are many investigations of the possible size of the
coefficients of a null-series. They show that the coefficients may
be arbitrarily close to $l^{2}$. See \cite{I57}, \cite{A84}, \cite{P85}, 
\cite{K87}. In
all known constructions the behavior of the amplitudes in the positive
and the negative parts of the spectrum is the same. \cite{KO03} shows
that a substantial nonsymmetry may occur. How far may this nonsymmetry
go? Is it possible for the anti-analytic amplitudes to decrease fast?
Equivalently, may a function in $\PLA\setminus H^{2}$ be smooth? 

The method used in \cite{KO03} is too coarse to approach this problem.
However, we proved recently that smooth and even $C^{\infty}$ functions
do exist in $\PLA\setminus H^{2}$. Precisely, in \cite{KO04} we
sketched the proof of the following:

\begin{thm}
\label{thm:cr}There exists a null-series {\rm (\ref{eq:Null})} with amplitudes
in negative spectrum {\rm (}$n<0${\rm )} satisfying the condition\begin{equation}
c(n)=O(|n|^{-k}),\quad k=1,2,\dotsc\, .\label{eq:Cinfty}\end{equation}
\end{thm}
\vskip8pt

Hence we are lead to the following question: \emph{what is the maximal
possible smoothness of a \/{\rm ``}\/nonclassic\/{\rm '' PLA}\/ function}? In other
words we want to characterize the possible rate of decreasing the
amplitudes $|c(n)|$ of a null-series as $n\to-\infty$. This is the
main problem considered here.

\demo{{\rm 1.3}} It should be mentioned that if one replaces convergence a.e.~by convergence
on a set of positive measure, then the characterization is given by
the classic quasi-analyticity condition. Namely, the class of series
(\ref{eq:Null}) satisfying \begin{equation}
c(n)=O(e^{-\rho(|n|)})\quad\forall n<0\label{eq:cnoewn}\end{equation}
for some $\rho(n)$ (with some regularity) is prohibited from containing
a nontrivial series converging to zero on a set $E$ of positive
measure if and only if\begin{equation}
\sum\frac{\rho(n)}{n^{2}}=\infty.\label{eq:clasquasi}\end{equation}
The part ``only if'' is well known: if this sum converges one
may construct a function vanishing on an interval $E$ whose Fourier
coefficients satisfy (\ref{eq:cnoewn}), and for $n$ positive as
well (see e.g.~\cite[Chap.\ 6]{M35}). The ``if'' part follows
from a deep theorem of Beurling \cite{B89}, extended by Borichev
\cite{B88}. See more details below in Section \ref{sub:uniqhist}.

It turns out that in our situation the threshold is completely different.
The following uniqueness theorem with a much weaker requirement on coefficients
is true.

\begin{thm}
\label{thm:unique}
Let $\omega$ be a function $\mathbb{R}^{+}\to\mathbb{R}^{+}${\rm ,}
$\omega(t)/t$ increase and\begin{equation}
\sum\frac{1}{\omega(n)}<\infty.\label{eq:sumwfin}\end{equation}
Then the condition\/{\rm :}
\begin{equation}
c(n)=O(e^{-\omega(\log|n|)}),\; n<0\label{eq:qa}\end{equation}
for a series {\rm (\ref{eq:Null})} converging to zero a.e.~implies that
all $c(n)$ are zero.
\end{thm}

It is remarkable that the condition is sharp. The following strengthened
version of Theorem \ref{thm:cr} is true:

\begin{thm}
\label{thm:precise}
Let $\omega$ be a function $\mathbb{R}^{+}\to\mathbb{R}^{+}${\rm ,} let
$\omega(t)/t$ be concave and\begin{equation}
\sum\frac{1}{\omega(n)}=\infty.\label{eq:sumwinf}\end{equation}
Then there exists a null-series {\rm (\ref{eq:Null})} such that {\rm (\ref{eq:qa})}
is fulfilled.
\end{thm}

So the maximal possible smoothness of a ``nonclassical'' $\PLA$
function $f$ is precisely characterized in terms of its Fourier transform
by the condition\[
\widehat{f}(n)=O(e^{-\omega(\log|n|)}), \quad n\in\mathbb{Z}\]
where $\omega$ satisfies (\ref{eq:sumwinf}). As far as we are aware
this condition has never appeared before as a smoothness threshold.

We mention that whereas the usual quasi-analyticity is placed near
the ``right end'' in the scale of smoothness connecting $C^{\infty}$
and analyticity, this new quasi-analyticity threshold is located just
in the opposite side, somewhere between $n^{-\log\log n}$ and $n^{-(\log\log n)^{1+\varepsilon}}$.

The main results of this paper were announced in our recent note \cite{KO04}.

\section{Preliminaries}

In this section we give standard notation, needed background and
some additional comments.

\demo{{\rm 2.1}} We denote by $\mathbb{T}$ the circle group $\mathbb{R}/2\pi\mathbb{Z}$.
We denote by $\mathbb{D}$ the disk in the complex plane $\{ z:|z|<1\}$
and $\partial\mathbb{D}=\{ e^{it}:t\in\mathbb{T}\}$. For a function
$F$ (harmonic, analytic) on $\mathbb{D}$ and a  $\zeta\in\partial\mathbb{D}$
we shall denote the nontangential limit of $F$ at $\zeta$ (if it
exists) by $F(\zeta)$.

We denote by $C$ and $c$ constants, possibly different in different
places. By $X\approx Y$ we mean $cX\leq Y\leq CX$. By $X\ll Y$
we mean $X=o(Y)$. Sometimes we will use    notation such as $-O(\cdot)$.
While this seems identical to just $O(\cdot)$ we use this notation
to remind the reader that the relevant quantity is negative.

The notation $\left\lfloor x\right\rfloor $ will stand for the lower
integral value of $x$. $\left\lceil x\right\rceil $ will stand for
the upper integral value.

When $x$ is a point and $K$ some set in $\mathbb{T}$ or $\mathbb{D}$,
the notation $d(x,K)$ stands, as usual, for $\inf_{y\in K}d(x,y)$.

\demo{{\rm 2.2}} For a $z\in\mathbb{D}$ we shall denote the Poisson kernel at the
point $z$ by $P_{z}$ and the conjugate Poisson kernel by $Q_{z}$.
We denote by $H$ the Hilbert kernel on~$\mathbb{T}$. See e.g.~\cite{Z59}.
If $f\in L^{2}(\mathbb{T})$ we shall denote by $F(z)$ the harmonic
extension of $f$ to the disk, i.e. \begin{equation}
F(z)=\int_{0}^{2\pi}P_{z}(t)f(t)\, dt,\quad\forall z\in\mathbb{D}.\label{eq:defPz}\end{equation}
Similarly, the harmonic conjugate to $F$ can be derived directly
from $f$ by \[
\widetilde{F}(z)=\int_{0}^{2\pi}Q_{z}(t)f(t)\, dt,\quad\forall z\in\mathbb{D}.\]
It is well known that $F$ and $\widetilde{F}$ have nontangential
boundary values a.e.~and that $F(e^{it})=f(t)$ a.e. We shall denote
$\widetilde{f}(t):=\widetilde{F}(e^{it})$. We remind the reader also that
\[
\widetilde{f}(x)=(f*H)(x)=\int_{0}^{2\pi}f(t)H(x-t)\, dt\]
where the integral is understood in the principal value sense.

For a function $F$ on the disk, the notation $F^{(D)}$ denotes tangent
differentiation, namely $F'(re^{i\theta}):=\frac{\partial F}{\partial\theta}$.
 \pagebreak
The representations above admit differentiation. For example,\[
F^{(D)}(z)=\int P_{z}^{(D)}(t)f(t)\, dt,\quad\forall z\in\mathbb{D}.\]
 We shall use the following well known estimates for $P$, $Q$ and
their derivatives: \begin{equation}
|P_{z}^{(D)}(t)|\leq\frac{C(D)}{|e^{it}-z|^{D+1}},\quad|Q_{z}^{(D)}(t)|\leq\frac{C(D)}{|e^{it}-z|^{D+1}}\quad\forall D\geq0;\label{eq:IPtag}\end{equation}
for $H$ we shall need the symmetry $H(t)=-H(-t)$ and \begin{equation}
|H^{(D)}(t)|\leq\frac{(CD)^{CD}}{|e^{it}-1|^{D+1}}.\label{eq:Hilbert}\end{equation}

\setcounter{Subsec}{2}
\Subsec{\label{sub:uniqhist}Uniqueness theorems}
In 1918 Privalov proved the following fundamental theorem:

\medbreak {\it Let $F$ be an analytic function on $\mathbb{D}$
such that $F(e^{it})=0$ on a set $E$ of positive measure. Then $F$
is identically zero.}
\medskip{}

See \cite{P50}, \cite{K98}. The conclusion also holds under the condition
$$
F(e^{it})=\sum_{n=-\infty}^{-1}c(n)e^{int}\quad\textrm{on }E$$
with the $|c(n)|$ decreasing exponentially. When one goes further
the picture gets more complicated. Examine the following result of
Levinson and Cartwright \cite{L40}:

\medskip{}
  \emph{Let $F$ be an analytic function on $\mathbb{D}$
with the growth condition\begin{equation}
|F(z)|<\nu(1-|z|)\quad\int_{0}^{1}\log\log\nu<\infty.\label{eq:loglog}\end{equation}
Assume that $F$ can be continued analytically through an arc $E\subset\partial\mathbb{D}$
to an $f$ in $\mathbb{C}\setminus\mathbb{D}$ which satisfies\[
f(z)=\sum_{n=-\infty}^{-1}c(n)z^{n}.\]
and the $c(n)$ satisfy the quasi-analyticity conditions {\rm (\ref{eq:cnoewn}),
(\ref{eq:clasquasi})}. Then $F$ and $f$ are identically zero.}
\medskip{}

It follows if a series (\ref{eq:Null}) converges to zero on an interval
and the ``negative'' coefficients decrease quasianalytically then
it is trivial.

In 1961 Beurling extended the Levinson-Cartwright theorem from an
arc to any set $E$ with positive measure (see \cite{B89}):

\medskip{}
 \emph{Let $f\in L^{2}$ vanish on $E$ and let its Fourier
coefficients $c(n)$ satisfy {\rm (\ref{eq:cnoewn}),~(\ref{eq:clasquasi})}.
Then $f$ is identically zero.}
\medskip{}

Borichev \cite{B88} proved that the $L^{2}$ condition in this theorem
could be replaced by a very weak growth condition on the analytic
part $F$ in $\mathbb{D}$, similar in spirit to (\ref{eq:loglog}).
Certainly this condition would be fulfilled if the series converged
pointwise on $E$. Note again that the classic quasi-analyticity condition
in all these results cannot be improved. Our proof of uniqueness
uses the same general framework used in \cite{B88}, \cite{BV89}, \cite{Bo89}.

Other results about the uncertainty principle in analytic settings
exist, namely connecting smallness of support with fast decrease of
the Fourier coefficients. See for example \cite{H78} for an analysis
of support of measures with smooth Cauchy transform. The connection
between the smoothness of the boundary value of a function $F$ and
the increase of $F$ near the singular points of the boundary was
investigated for $F$ from the Nevanlinna class; see Shapiro \cite{S66},
Shamoyan \cite{S95} and Bourhim, El-Fallah and  Kellay \cite{BEK04}.
In particular, applying theorem A of \cite{BEK04} to our case shows
that one cannot construct a $C^{1}$ function in $\PLA\setminus H^{2}$
by taking the boundary value of a Nevanlinna function. For comparison,
our first example of a function from $\PLA\setminus H^{2}$ (see \cite{KO03})
is a boundary value of a Nevanlinna class function. That example is
$L^{\infty}$ and can be made continuous, but it cannot be made smooth
in any reasonable sense without leaving the Nevanlinna class.

\Subsec{The harmonic measure}
Let $\mathcal{D}$ be a connected open set in $\mathbb{C}$ such that
$\partial\mathcal{D}$ is a finite collection of Jordan curves, and
let $v\in\mathcal{D}$. Let $B$ be Brownian motion (see \cite[I.2]{B95})
starting from $v$. Let $T$ be the stopping time on the boundary
of $\mathcal{D}$, i.e.\[
T:=\inf\{ t:B(t)\in\partial\mathcal{D}\}.\]
See \cite[Prop.\ I.2.7]{B95}. Then $B(T)$ is a random point
on $\partial\mathcal{D}$, or in other words, the distribution of
$B(T)$ is a measure on $\partial\mathcal{D}$ called the \emph{harmonic
measure} and denoted by $\Omega(v,\mathcal{D})$. The following result
is due to Kakutani \cite{K44}.

\medskip{}
  \emph{Let $f$ be a harmonic function in a domain $\mathcal{D}$
and continuous up to the boundary. Let $v\in\mathcal{D}$. Then\begin{equation}
f(v)=\int f(\theta)\, d\Omega(v,\mathcal{D})(\theta).\label{eq:kak}\end{equation}
}

It follows that the definition of harmonic measure above is equivalent
to the original definition of Nevanlinna which used solutions of Dirichlet's
problem. We shall also need the following version of Kakutani's theorem:

\medskip{}
 \emph{Let $f$ be a subharmonic function in a domain $\mathcal{D}$
and upper semi-continuous up to $\partial\mathcal{D}$. Let $v\in\mathcal{D}$.
Then}
\begin{equation}
f(v)\leq\int f(\theta)\, d\Omega(v,\mathcal{D})(\theta).\label{eq:subkak}\end{equation}

See \cite[Propositions II.6.5 and II.6.7]{B95}. See also   [{\it ibid},
Theorem II.1.15 and   Proposition II.1.13].

\section{\label{sec:Construction}Construction of smooth $\PLA$ functions}

 3.1.\quad  In this section we prove Theorem \ref{thm:precise}. We wish to restate
it in a form which makes explicit the fact that the singular set is
in fact compact:

\demo{\scshape Theorem \ref{thm:precise}$'$}
{\it Let $\omega$ be a function $\mathbb{R}^{+}\to\mathbb{R}^{+}${\rm ,} $\omega(t)/t$ be
concave and $\sum\frac{1}{\omega(n)}=\infty$. Then there exists a
series {\rm (\ref{eq:Null})} converging to zero outside a compact set $K$
of measure zero such that {\rm (\ref{eq:qa})} is fulfilled.}
\Enddemo

The regularity condition that $\omega(t)/t$ be concave in Theorem
\ref{thm:precise}$'$ implies the very rough estimate $\omega(t)=e^{o(t)}$,
which is what we will use. Actually, one may streng\-then the theorem
slightly by requiring only that $\omega(t)/t$ is increasing and $\omega(t)=e^{o(t)}$,
and the result would still hold.

Without loss of generality it is enough to prove \begin{equation}
c(n)=O(e^{-c\omega(\log|n|)}),\; n<0\label{eq:qaC}\end{equation}
 for some $c>0$, instead of (\ref{eq:qa}). Also we may assume $\omega(t)/t$
increases to infinity (otherwise, just consider $\omega(t)=t\log(t+2)$
instead).

The $c$ above, like all notation $c$ and $C$, $\ll$, $o$ and
$O$ in this section, is allowed to depend on $\omega$. In general
we will consider $\omega$ as given and fixed, and will not remind
the reader that the various parameters depend on it.

A rough outline of the proof is as follows: we shall define a probabilistically-skewed
thick Cantor set $K$ and a random harmonic function $G$ on the disk
such that the boundary values of $G$ on $K$ are positive infinite,
while the boundary values outside $K$ are finite negative (except
a countable set of points where they are infinite negative). Further,
the function $G$ is ``not integrable'' in the sense that $\int_{0}^{2\pi}|G(re^{i\theta})|\, d\theta\rightarrow\infty$
as $r\to1$. The thickness of the set $K$ would depend on $\omega$.
For example, if $\omega(t)=t\log t$ (which is enough for the construction
of a nonclassic $\PLA\cap\, C^{\infty}$ function, i.e.~for the proof
of Theorem \ref{thm:cr}) then $K$ would have infinite $\delta\log\log1/\delta$-Hausdorff
measure. Then we shall define $F=e^{G+i\widetilde{G}}$ and $f$ its
boundary value ($f$ is a nonclassic $\PLA$ function). We shall
arrange for $G|_{\partial\mathbb{D}}$ to converge to $-\infty$ sufficiently
fast near $K$, and it would follow that $f$ is smooth. A bound for
the growth of $G$ to $+\infty$ near $K$ would ensure that the Taylor
coefficients of $F$ go to zero with probability one. Finally the
desired null-series would be defined by \begin{equation}
c(n):=\widehat{f}(n)-\begin{cases}
\widehat{F}(n) & n\geq0\\
0 & n<0\end{cases}\label{eq:defcnfF}\end{equation}
where $\widehat{f}$ is the Fourier transform of $f$ while $\widehat{F}$
are the Taylor coefficients of $F$:\begin{equation}
F(z)=\sum_{n=0}^{\infty}\widehat{F}(n)z^{n}.\label{eq:Taylor}\end{equation}

\setcounter{Subsec}{1}
\Subsec{\noindent Auxiliary sequences}
 Let $\omega_{2}$ satisfy that $\omega_{2}(t)/t$ is increasing,
$\sum\frac{1}{\omega_{2}}\break =\infty$ , and \begin{equation}
\omega(t)\ll\omega_{2}(t)=\omega(t)t^{o(1)}\label{eq:defmun}\end{equation}
(note that  $\omega_{2}(t)=e^{o(t)}$). Define \[
\Phi(n):=\exp\left(-\sum_{k=1}^{n}\frac{1}{\omega_{2}(k)}\right) \]
and in particular $\Phi(0)=1$. Also,  $\Phi$ decreases slowly (depending on
$\omega_{2}$), and  the fact that $\frac{\omega_{2}(t)}{t}$
increases to $\infty$ gives \begin{equation}
\Phi(n)=n^{-o(1)}.\label{eq:Phinsmall}\end{equation}
Another regularity condition over $\Phi$ that will be used is the
following:

\begin{lem}
\label{lem:Phinormal}\begin{equation}
\sum_{k=1}^{n}\Phi(k)=O(n\Phi(n)).\label{eq:sumPhiPhi}\end{equation}
\end{lem}

\Proof 
Fix $N$ such that $\omega_{2}(n)/n>100$ for $n>N$ . Then for all
$n>3N$,\[
\sum_{k=\left\lceil \frac{1}{3}n\right\rceil }^{n}\frac{1}{\omega_{2}(k)}\leq0.03\]
and hence $\Phi(k)\leq1.04\Phi(n)$ for all $k\in\left[\left\lfloor \frac{1}{3}n\right\rfloor ,n\right]$.
Inductively we get $\Phi(k)\leq1.04^{l}\Phi(n)$ for any $k\in\left[\left\lfloor n3^{-l}\right\rfloor ,\left\lfloor n3^{1-l}\right\rfloor \right]\cap\{ N,\dotsc\}$.
Hence we get\begin{align*}
\sum_{k=1}^{n}\Phi(k) & =\sum_{l=1}^{\left\lfloor \log_{3}n\right\rfloor +1}\sum_{k=\left\lfloor n3^{-l}\right\rfloor +1}^{\left\lfloor n3^{1-l}\right\rfloor }\Phi(k)\leq N+\sum_{l=1}^{\left\lfloor \log_{3}n\right\rfloor +1}\left\lceil 2\cdot3^{-l}n\right\rceil \Phi(n)(1.04)^{l}\\
 & =O(n\Phi(n)).\end{align*}
Notice that in the last equality we used the fact that $\Phi(n)\gg1/n$ (\ref{eq:Phinsmall}).
\Endproof\vskip4pt  

 Next, define
 \begin{equation}
\sigma_{n}=2\pi\cdot2^{-n}\Phi(n),\quad\tau_{n}=\frac{1}{12}(\sigma_{n-1}-2\sigma_{n}),
\quad
n\geq0.\label{eq:deftaun}\end{equation}
 The purpose behind the definition of $\Phi$ is so that the following
(which can be verified with a simple calculation) holds: \begin{equation}
\frac{\tau_{n}}{\sigma_{n}}=\frac{1}{6\omega_{2}(n)}+O\left(\frac{1}{\omega_{2}^{2}(n)}\right).\label{eq:tausigmanu3}\end{equation}
From this and the regularity conditions $\omega_{2}(n)=e^{o(n)}$
and (\ref{eq:Phinsmall}) we get a rough but important estimate for
$\tau_{n}$: \begin{equation}
\tau_{n}=2^{-n-o(n)}.\label{eq:taunreg}\end{equation}

\Subsec{The functions $g_{n}$}
Next we define some auxiliary functions. Let $a\in C^{\infty}([0,1])$
be a nonnegative function satisfying \[
a|_{\left[0,1/3\right]}\equiv0,\quad a|_{[1/2,1]}\equiv1,\quad\max\left|a^{(D)}\right|\leq(CD)^{CD}.\]
Since the standard building block $e^{-1/x}$ satisfies the estimate
for the growth of the derivatives above (even a very rough estimate
can show this --- say, use Lemma \ref{lem:expderivk} below), and since
such constraints are preserved by multiplication, there is no difficulty
in constructing $a$. 

Let $l$ be defined by \[
l(t)=-t^{-1/3}a(1-t)-a(t).\]
Then $l$ satisfies \begin{alignat}{2}
l(t) & =-t^{-1/3},\qquad & t & \in\left]0,{\textstyle \frac{1}{3}}\right],\label{eq:lsing}\\
l(t) & =-1, & t & \in[{\textstyle \frac{2}{3}},1],\nonumber \end{alignat}
and $l\leq-1$ on $\left]0,1\right]$.

Using $l$, define functions on $\mathbb{R}$ depending on a parameter
$s\in[0,1]$, \begin{equation}
l^{\pm}(s;x):=\begin{cases}
l(x) & 0<x\leq1\\
-1 & 1<x\leq2\pm s\\
l(3\pm s-x) & 2\pm s<x\leq3\pm s\end{cases}\label{eq:deflpm}\end{equation}
and $0$ otherwise. The estimate for the derivatives of $a$ translates
to \begin{equation}
\left|\left(l^{\pm}\right)^{(D)}(s;x)\right|\leq\frac{(CD)^{CD}}{d(x,\{3\pm s,0\})^{D+1/3}}.\label{eq:lverysmth}\end{equation}

Let $s(n,k)$ be a collection of numbers between $0$ and $1$, for
each $n\in\mathbb{N}$ and each $0\leq k<2^{n}$. Most of the proof
will hold for any choice of $s(n,k)$, but in the last part we shall
make them random, and prove that the constructed function will have
the required properties for almost any choice of $s(n,k)$. Define
now inductively intervals $I(n,k)=[a(n,k),a(n,k)+\sigma_{n}]$ (we
call these $I(n,k)$ ``intervals of rank $n$'') as follows: $I(0,0)=[0,2\pi]$
and for $n\geq0$, $0\leq k<2^{n}$,
 \begin{align}
a(n+1,2k) & =a(n,k)+\tau_{n+1}(3+s(n+1,2k)),\label{eq:defank} \\
a(n+1,2k+1) & =a(n,k)+{\textstyle 
\frac{1}{2}}\sigma_{n}+\tau_{n+1}(3+s(n+1,2k+1)).\nonumber\end{align}
In other words, at the $n^{\rm th}$ step, inside each interval of rank
$n$ (which has length $\sigma_{n}$), situate two disjoint intervals
of rank $n+1$ of lengths $\sigma_{n+1}$ in random places (but not
too near the boundary of $I(n,k)$ or its middle). Define \begin{align*}
K & :=\cap_{n=1}^{\infty}K_{n}, & K_{n} & :=\cup_{k=0}^{2^{n}-1}I(n,k).\\
K^{\circ} & :=e^{iK}, & K_{n}^{\circ} & :=e^{iK_{n}}.\end{align*}
Note that $\sum\frac{1}{\omega_{2}}=\infty$ shows that $\Phi(n)\to0$
and hence $K$ has zero measure.

We now define the most important auxiliary functions, $g_{n}\in L^{2}(\mathbb{T})$.
We define them inductively, with $g_{0}\equiv0$. For one $n$ and
$k$, let $I$ be the interval of rank $n-1$ containing $I(n,k)$,
and let $I'$ be its half containing $I(n,k)$. Now, $I'\setminus I(n,k)$
is composed of two intervals, which we denote by $J_{1}$ (left) and
$J_{2}$ (right). Define the function $g_{n}(t)$ on the set $I'\setminus I(n,k)$
by
\begin{equation}
\begin{aligned}g_{n}(t) & :=\omega(n)l^{+}(s(n,k);\varphi_{1}(t)), \enspace & t & \in J_{1}, \enspace& 
&
\varphi_{1}:J_{1}\to[0,3+s(n,k)],\\ g_{n}(t) & :=\omega(n)l^{-}(s(n,k);\varphi_{2}(t)),\enspace & t & \in
J_{2}, \enspace&  &
\varphi_{2}:J_{2}\to[0,3-s(n,k)],\end{aligned}
\label{eq:defgn}\end{equation}
where the $\varphi_{i}$-s are linear, increasing and onto so that they
are defined uniquely by their domain and range. As will become clear
later, the $\omega(n)$ factor above is what determines the rate of
decrease of the coefficients of the null series.

\begin{figure}
\centerline{\epsfig{file=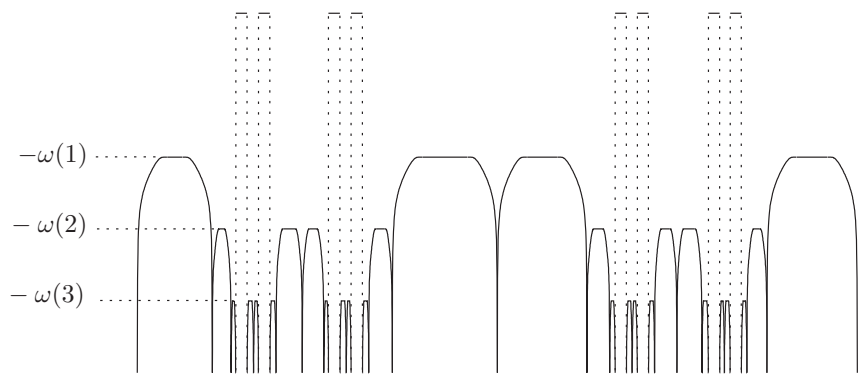}}
\caption{\label{cap:g}$g_{3}$ (not drawn to scale). Notice the random perturbations
in the widths of the constant parts of $g_{3}^{-}$ but the fixed
width of the intervals in $K_{3}$.}
\end{figure}

This defines $g_{n}$ on $K_{n-1}\setminus K_{n}$. On $\mathbb{T}\setminus K_{n-1}$
we define $g_{n}\equiv g_{n-1}$. On $K_{n}$ we define $g_{n}$ to
be a constant such that $\int_{\mathbb{T}}g_{n}=0$. Note that $g_{n}$
is negative on $\mathbb{T}\setminus K_{n}$ and positive on $K_{n}$.
Also note that the definition of $g_{n}$ shows that $\int_{I(n-1,k)}g_{n}^{-}$
is independent of $s(n-1,k)$ --- what you earn on the left you lose
on the right. See Figure \ref{cap:g}.%

Extend $g_{n}(e^{it})$ to a harmonic function in the interior of
the disk (remember that each $g_{n}$ is in $L^{2}$), and denote
the extension by $G_{n}$. Denote by $\widetilde{G_{n}}$ the harmonic
conjugate to $G_{n}$.

\Subsec{The growth of the $g_{n}$}
We need to estimate the positive part of $g_{n}$. We have\begin{equation}
\int_{I(n-1,k)}|g_{n}^{-}(x)|\stackrel{(*)}{\approx}
\tau_{n}\omega(n)\stackrel{(**)}{\ll}
\sigma_{n}\label{eq:vkapprox}\end{equation}
where $(*)$ comes from the definition of $g_{n}$ (\ref{eq:defgn})
and $(**)$ comes from $\tau_{n}/\sigma_{n}\approx1/\omega_{2}(n)$
(\ref{eq:tausigmanu3}) and $\omega\ll\omega_{2}$ (\ref{eq:defmun}).
Summing (and using $\sigma_{n}=2^{-n}\Phi(n)$) we get \begin{equation}
\int_{\mathbb{T}}\left|g_{n}^{-}\right|=\sum_{l=0}^{n-1}\sum_{k=0}^{2^{l}-1}
\int_{I(l,k)}|g_{l+1}^{-}(x)|\stackrel{(\ref{eq:vkapprox})}{=} 
\sum_{l=1}^{n}o(\Phi(l))\stackrel{(*)}{=}o(n\Phi(n)),
\label{eq:Wnapprox}\end{equation}
where $(*)$ comes from Lemma \ref{lem:Phinormal} and $\Phi(n)\gg1/n$
(\ref{eq:Phinsmall}). Hence \begin{equation}
\max g_{n}=o(n).\label{eq:maxfn}\end{equation}
This crucial inequality is the one that guarantees in the end that
our function $F$ would satisfy $\widehat{F}(m)\to0$. Comparing this  to
(\ref{eq:defgn}) we observe that even though $K$ has zero measure,
one can balance superlinear growth outside $K$ (the $\omega(n)$
factor in (\ref{eq:defgn})) with sublinear growth inside $K$.

We will also need a simple estimate from the other side. The same
calculations, but using $\omega(n)=\omega_{2}(n)\cdot n^{-o(1)}$
(the second half of (\ref{eq:defmun})) and $\Phi(n)=n^{-o(1)}$ (\ref{eq:Phinsmall})
give\begin{equation}
\int_{\mathbb{T}}g_{n}^{-}=-n^{1-o(1)}.\label{eq:l1normfn}\end{equation}

\Subsec{The limit of the $G_{n}$}
First we want to show that the $G_{n}$'s converge to a
harmonic function $G$ on compact subsets of the disk, and to discuss
the boundary behavior of $G$ and $\widetilde{G}$. For this purpose
we need to examine the singularities of $g_{n}$. First, and more
important is $K$. Clearly, $\lim_{n\rightarrow\infty}g_{n}(t)=+\infty$
while $\lim_{t'\rightarrow t,t'\not\in K} \lim_{n\rightarrow\infty}g_{n}(t')=-\infty$
for every $t\in K$. Additionally we have a countable set of points
where the $g_{n}$'s have $t^{-1/3}$-type singularities, namely \[
Q:=\bigcup_{n=1}^{\infty}\bigcup_{k=0}^{2^{n}-1}\left\{ a(n,k),a(n,k)+{\textstyle \frac{1}{2}}\sigma_{n},a(n,k)+\sigma_{n}\right\} .\]
Denote $K':=K\cup Q$, $(K')^{\circ}:=e^{iK'}$. We will need the
following calculation:

\begin{lem}
For any $z\in\overline{\mathbb{D}}\setminus K_{n}^{\circ}${\rm ,} and any
$D\geq0${\rm ,}
\begin{equation}
|G_{n+1}^{(D)}(z)-G_{n}^{(D)}(z)|\leq\frac{C(D)}{2^{n}d(z,K_{n}^{\circ})^{D+1}}.\label{eq:fnDfn1DdzK}\end{equation}
\end{lem}

\Proof 
On the circle $\mathbb{T}$, $g_{n+1}-g_{n}$ is nonzero only on
the intervals $I(n,k)$, and on each interval we have\begin{equation}
\int_{I(n,k)}\left(g_{n+1}(x)-g_{n}(x)\right)\, dx=0.\label{eq:fnfnmin1zero}\end{equation}
Further, the negative part of $g_{n+1}-g_{n}$ on $I(n,k)$, which
is simply $g_{n+1}^{-}-\max g_{n}$ restricted to $I(n,k)\setminus(I(n+1,2k)\cup I(n+1,2k+1))$
can be estimated using (\ref{eq:vkapprox}) and (\ref{eq:maxfn})
to get \begin{align*}
\int_{I(n,k)}|g_{n+1}(x)-g_{n}(x)| & =2\int(g_{n+1}-g_{n})^{-}\approx\tau_{n+1}(\omega(n+1)+o(n))\\
 & \stackrel{(*)}{\ll}2^{-n}\Phi(n+1)\ll2^{-n}\end{align*}
where in $(*)$ we use that $n\ll\omega(n)$. Hence by (\ref{eq:fnfnmin1zero}),\begin{equation}
\Big|\int_{t}^{u}g_{n+1}(x)-g_{n}(x)\, dx\Big|\leq C2^{-n}\quad\forall t,u\in[0,1],\,\forall n.\label{eq:intgparts}\end{equation}
Write $G^{(D)}(z)=\int_{\mathbb{T}}g(t)P_{z}^{(D)}(t)$, where $P_{z}$
is the Poisson kernel. We divide into two cases: if $1-|z|>\frac{1}{2}d(z,K_{n}^{\circ})$
then we have from (\ref{eq:IPtag}) that\[
\int_{\mathbb{T}}|P_{z}^{(D+1)}|\leq\frac{C(D)}{(1-|z|)^{D+1}}\leq\frac{C(D)}{d(z,K_{n}^{\circ})^{D+1}}.\]
On the other hand, if $1-|z|\leq\frac{1}{2}d(z,K_{n}^{\circ})$ then
$g_{n+1}-g_{n}$ is zero in an interval $J:=\left[t-cd(z,K^{\circ}),t+cd(z,K^{\circ})\right]$
for some $c$ sufficiently small, where $t$ is given by $e^{it}=z/|z|$,
and \[
\int_{\mathbb{T}\setminus J}|P_{z}^{(D+1)}|\leq\frac{C(D)}{d(z,K_{n}^{\circ})^{D+1}}.\]
In either case, a simple integration by parts gives (\ref{eq:fnDfn1DdzK}).

Finally, on $\partial\mathbb{D}$ we have $G_{n+1}(e^{it})-G_{n}(e^{it})=g_{n+1}(t)-g_{n}(t)=0$
for every $t\not\in K_{n}$.
\Endproof\vskip4pt  
A similar calculation with the conjugate Poisson kernel (and the Hilbert
kernel on the boundary) shows\begin{equation}
|\widetilde{G_{n+1}}^{(D)}(z)-\widetilde{G_{n}}^{(D)}(z)|\leq\frac{C(D)}{2^{n}d(z,K_{n}^{\circ})^{D+1}}.\label{eq:zntilde}\end{equation}
From (\ref{eq:fnDfn1DdzK}) and (\ref{eq:zntilde}) it is now clear
that both $G_{n}$ and $\widetilde{G_{n}}$ converge uniformly on
compact subsets of $\overline{\mathbb{D}}\setminus(K')^{\circ}$.
Denote their respective limits by $G$ and $\widetilde{G}$ --- clearly
they are indeed harmonic conjugates which justifies the notation $\widetilde{G}$.
Also we remind the reader the known fact that if $g_{n}$ is $C^{D}$
in some interval $I\subset\mathbb{T}$ then $G_{n}$ is $C^{D}$ in
$e^{iI}$ and in particular $G_{n}^{(D)}$ is continuous there. The
following lemma is now clear:

\begin{lem}
 \label{lem:cont}
\begin{itemize}
\ritem{(i)}   $G+i\widetilde{G}$ is analytic in $\mathbb{D}$
and continuous up to the boundary except at $(K')^{\circ}$.
\ritem{(ii)}  If $t\in\mathbb{T}\setminus K_{n}$ then $G(e^{it})=g_{n}(t)=g_{n}^{-}(t)$.
\ritem{(iii)}  $(G_{n}+i\widetilde{G_{n}})^{(D)}$ converges
to $(G+i\widetilde{G})^{(D)}$ uniformly on compact subsets of $\overline{\mathbb{D}}\setminus(K')^{\circ}$.
\end{itemize}
\end{lem}

\vglue-12pt
\Subsec{\noindent The function $F$}
We can now define a crucial element of the construction:
\[
F=\exp(G+i\widetilde{G}).\]
Clearly Lemma \ref{lem:cont}, (i) shows that $F$
is an analytic function with almost everywhere defined boundary values.
Denote by $f(t)$ the boundary value of $F$ at $e^{it}$. Define
similarly $g$ and $\widetilde{g}$ and get that $f=e^{g+i\widetilde{g}}$.

The reader should keep in mind that the relation between $F$ and
$f$ is not similar to the one between an $H^{2}$ function and its
boundary value (for example, between $G_{n}+i\widetilde{G_{n}}$ and
$g_{n}+i\widetilde{g_{n}}$). In our case there is a singular distribution
(supported on $K'$) which is ``lost'' when taking the  limit. The
Fourier series of this singular distribution is exactly the null series
we are trying to construct.

\begin{lem}
\label{lem:PLA} 
\begin{itemize}
\ritem{(i)} $F$ is not in $H^{1}(\mathbb{D})$.
\ritem{(ii)} $f\in L^{\infty}(\mathbb{T})$.
\end{itemize}
\end{lem}

The first follows from Lemma \ref{lem:cont}, (ii) if we
notice that the $L^{1}$ norms of $g_{n}$ tend to $\infty$ according
to (\ref{eq:l1normfn}) so that  $\log|f|=g\not\in L^{1}(\mathbb{T})$.
The second is also a direct consequence of Lemma \ref{lem:cont},
(ii). These properties taken together show that the $c(n)$
(\ref{eq:defcnfF}) are nontrivial. The theorem now divides into
the following two claims:

\begin{lem}
\label{lem:Fstarsmooth}$f\in C^{\infty}(\mathbb{T})$. Further{\rm ,} $f$
has the smoothness in the statement of the theorem\/{\rm :}
\begin{equation}
\widehat{f}(n)=O(e^{-c\omega(\log|n|)}),\; n\in\mathbb{Z}.\label{eq:Fstarsmooth}\end{equation}
\end{lem}

\begin{lem}
\label{lem:Taylor0}With probability $1${\rm ,}\[
\widehat{F}(n)=o(1).\]
\end{lem}

 Now we use the Riemann localization principle in a form due to Kahane-Salem
\cite[p.~54]{KS94}: 

\medskip{}
  \emph{If $S$ is a distribution with $\widehat{S}(n)=o(1)$
and $I$ is an interval outside the support of 
$S$ then $\sum\widehat{S}(n)e^{int}=0$
on $I$. }
\medskip{}

Lemma \ref{lem:cont}, (i) shows that the $c(n)$
defined by (\ref{eq:defcnfF}) represent a singular distribution supported
on $K'$, and the last estimate shows that $c(n)=o(1)$. Hence (\ref{eq:Null})
is a nontrivial series convergent to $0$ everywhere on $\mathbb{T}\setminus K'$.
This, along with Lemma \ref{lem:Fstarsmooth}, proves the theorem.

The purpose of the next section is to prove Lemma \ref{lem:Fstarsmooth}.

\Subsec{Smoothness}
The following two lemmas are self-contained; that is, their $f$-s,
$g$-s, $\omega$-s and $K$-s are not necessarily the same ones as
those defined in the previous parts of the proof.

\begin{lem}
\label{lem:expderivk}Let $f=\exp(g)$. Then \begin{equation}
f^{(D)}=f\cdot\sum_{l_{1}+\dotsc+l_{i}=D}a_{\vec{l}}\prod_{j=1}^{i}g^{(l_{j})},\label{eq:expderivk}\end{equation}
and $\sum_{\vec{l}}|a_{\vec{l}}|\leq D!$
\end{lem}

This is a straightforward induction and we shall skip the proof.

\begin{lem}
\label{lem:fsmthFsmth}Let $\omega(t)$ satisfy that $\omega(t)/t$
is increasing to $\infty$ and $\omega(t)=e^{o(t)}$. Let $K$ be
some compact and let $g\in C^{\infty}(\mathbb{T}\setminus K)$ satisfy
\begin{enumerate}
\item \label{enu:flow}$\real g(x)\leq-\omega(\log1/d(x,K))${\rm ;}
\item \label{enu:reg}$|g^{(D)}(x)|\leq\frac{(CD)^{CD}}{d(x,K)^{2D}}$ for
every $D\geq1$.
\end{enumerate}
Let $f=e^{g}$ outside $K${\rm ,} $f|_{K}\equiv0$. Then $\widehat{f}(m)=O(e^{-c\omega(\log|m|)})$.
\end{lem}

We remark that condition \ref{enu:reg} interfaces only with the regularity
condition $\omega=e^{o(t)}$. The important point here is the interaction
between condition \ref{enu:flow} and the estimate for $\widehat{f}$.

\Proof 
Denote $d=d(x,K)$. Plugging the inequality for $g^{(D)}$ into (\ref{eq:expderivk})
gives\[
|f^{(D)}(x)|\leq|f(x)|D!\frac{(CD)^{CD}}{d^{2D}}\leq|f(x)|\frac{(CD)^{CD}}{d^{2D}}.\]
In particular, $\omega(t)/t\to\infty$ shows that $|f(x)|\leq Ce^{-\omega(\log1/d)}$
decreases superpolynomially near $K$ which shows that $f^{(D)}(x)=0$
for all $x\in K$ and (inductively) for all $D$ and hence $f\in C^{\infty}([0,1])$.
Further (assume $D>1$),\begin{equation}
|f^{(D)}(x)|\leq C\exp\left(-\omega(\log1/d)+CD\log D+2D\log1/d\right).\label{eq:Fkx}\end{equation}
For any $m$ sufficiently large, choose now \[
D=\left\lfloor 2\frac{\omega(\frac{1}{4}\log|m|)}{\log|m|}\right\rfloor .
 \pagebreak
 \]
 Note that the condition $\omega(t)=e^{o(t)}$ shows that $D=|m|^{o(1)}$.
To estimate the maximum  of $f^{(D)}$ in (\ref{eq:Fkx}), we notice
that if $\log1/d>\frac{1}{4}\log|m|$ then $\omega(\log1/d)\geq2D\log1/d$
(here   $\omega(t)/t$ is increasing); hence we may estimate
roughly that \[
\max_{d}-\omega(\log1/d)+2D\log1/d\leq{\textstyle \frac{1}{2}}D\log|m|\]
and get\[
\Vert f^{(D)}\Vert _{\infty}\leq C\exp\big({\textstyle \frac{1}{2}}D\log|m|+CD\log D\big)\stackrel{(*)}{=}\exp\left(D\log|m|\left({\textstyle \frac{1}{2}}+o(1)\right)\right)\]
where $(*)$ comes from $D=|m|^{o(1)}$. 

We now use the fact that $|\widehat{f}(m)|\leq|m|^{-D}\Vert f^{(D)}\Vert _{\infty}$
to get\begin{align*}
|\widehat{f}(m)| & \leq C\exp(-({\textstyle \frac{1}{2}}-o(1))D\log|m|)\\
 & =C\exp(-\left(1-o(1)\right)\omega({\textstyle \frac{1}{4}}\log|m|)+O(\log|m|)).\end{align*}
Remembering that $\omega(\frac{1}{4}\log|m|)\leq\frac{1}{4}\omega(\log|m|)$
(again, because $\omega(t)/t$ is increasing) and that $\omega(t)/t\to\infty$ we see that
the lemma is proved.
\Endproof\vskip4pt  

We remark that, in some sense, the lemma actually hides two applications
of the Legendre transform, $(\mathfrak{L}h)(x):=\max_{t}h(t)-xt$.
Roughly speaking, the norms of $f^{(D)}$ are the Legendre transform
of the rate of decrease of $g$ to $-\infty$ (condition \ref{enu:flow}
of the lemma) and $\widehat{f}(m)$ are the Legendre transform of
$f^{(D)}$. Combining both facts allowed us not to use explicitly
the notation $\mathfrak{L}$ and  to simplify somewhat.

\demo{Proof of Lemma {\rm \ref{lem:Fstarsmooth}}} Our goal is to use Lemma \ref{lem:fsmthFsmth}
with the function $g+i\widetilde{g}$, the compact $K'$ and the $\omega$
of the lemma being $c\omega$ for some $c>0$ sufficiently small.
The condition \ref{enu:flow} on the size of the negative decrease
of $\real(g+i\widetilde{g})=g$ is easiest to show. Let $x\in K_{n-1}\setminus K_{n}$.
We divide into two cases: if $d(x,K')>e^{-n}$ then we may estimate\begin{equation}
-\omega(\log1/d(x,K'))\geq-\omega(n)\stackrel{(\ref{eq:defgn})}{\geq}g_{n}^{-}(x)\stackrel{(*)}{=}g(x)\label{eq:forlem1a}\end{equation}
where $(*)$ comes from Lemma \ref{lem:cont}, (ii). If
$d(x,K')\leq e^{-n}$ then $\tau_{n}=e^{-n(\log2+o(1))}\break \geq cd(x,K')^{0.7}$
(\ref{eq:taunreg}) and \begin{align}
g(x) & \leq-c\omega(n)\left(\frac{d(x,K')}{\tau_{n}}\right)^{-1/3}
\stackrel{(*)}{\leq}-cd(x,K')^{-0.1}\label{eq:forlem1} \\
 & =-c\exp\left({\textstyle \frac{1}{10}}\log1/d(x,K')\right)\stackrel{(**)}{\ll}-\omega(\log1/d(x,K'))
\nonumber\end{align}
where in $(*)$ we estimated trivially $\omega(n)\geq c$ and in $(**)$
we used the regularity condition $\omega(n)=e^{o(n)}$. Hence we get
$g(x)\leq-c\omega(\log1/d(x,K'))$ for all $x$, i.e.~the condition
\ref{enu:flow} of \pagebreak  Lemma \ref{lem:fsmthFsmth}.

To estimate $g^{(D)}$ outside $K'$, start from (\ref{eq:lverysmth})
and get for $x\in K_{n-1}\setminus K_{n}$ that \[
\left|g_{n}^{(D)}(x)\right|\leq\omega(n)\frac{(CD)^{CD}}{d(x,K')^{D+1/3}}.\]
Since $|I_{n,k}|\leq2\pi\cdot2^{-n}$ we get that for every $x\in I_{n,k}$
\[
d(x,K')^{2/3}\leq C2^{-2n/3}\ll1/\omega(n)\]
so that
\begin{equation}
\left|g^{(D)}(x)\right|=\left|g_{n}^{(D)}(x)\right|\leq\frac{(CD)^{CD}}{d(x,K')^{D+1}}.\label{eq:fstarderivd}\end{equation}
Note that (\ref{eq:fstarderivd}) holds for $g_{n}$ and any $x\not\in K_{n}$
(not necessarily in $K_{n-1}$).

For $\widetilde{g}$ we need to examine $\widetilde{g_{n}}$ and take
$n\to\infty$ (remember Lemma \ref{lem:cont}, (iii)).
Let $x\not\in K_{n}\cup K'$, let $\rho=\frac{1}{2}d(x,K')$ and let
$I=[-\rho,\rho]$. Now,
 \begin{equation}
\widetilde{g_{n}}(x)=\int_{\mathbb{T}}H(t)g_{n}(x-t)\, dt=\Big(\int_{I}+\int_{\mathbb{T}\setminus I}\Big)H(t)g_{n}(x-t)\, dt.\label{eq:intIint1mI}\end{equation}
In general, the symmetry of $H$ and $|H(t)|\approx\frac{1}{t}$ (\ref{eq:Hilbert})
allows to estimate for any~$h$
\begin{equation}
\left|\int_{I}h(t)H(t)\right|\leq C|I|\max_{I}|h'|\label{eq:HilbSob}\end{equation}
which we use as follows: For the first integral, $D$ differentiations
under the integral sign (which can be justified easily using (\ref{eq:HilbSob}))
show that\begin{align}
\left|\left(\frac{d}{dx}\right)^{D}\int_{I}H(t)g_{n}(x-t)\right| & 
=\left|\int_{I}H(t)\left(g_{n}\right)^{(D)}(x-t)\right|
\stackrel{(\ref{eq:HilbSob})}{\leq}C\rho\max_{I}\left|\left(g_{n}\right)^{(D+1)}\right|
\label{eq:intIf}\\
 & \stackrel{(\ref{eq:fstarderivd})}{\leq}\frac{(CD)^{CD}}{\rho^{D+1}}.
\nonumber\end{align}
For the second half of (\ref{eq:intIint1mI}) we consider $g_{n}$
and $H$ as periodic functions and change variables to get $\int_{x+\rho}^{x+2\pi-\rho}H(x-t)g_{n}(t)$.
This we differentiate $D$ times under the integral sign and shift
back, and we get
\begin{align} \label{eq:Int01notI}
\left(\frac{d}{dx}\right)^{D}\int_{\mathbb{T}\setminus I}H(t)g_{n}(x-t)  
=&\sum_{i=0}^{D-1}\left.H^{(i)}(t)\left(g_{n}\right)^{(D-1-i)}(x-t)\right|_{-\rho}^{\rho} \\
 &+\int_{\mathbb{T}\setminus I}H^{(D)}(t)g_{n}(x-t)\, dt\nonumber
\end{align}
where as usual $\left.g\right|_{a}^{b}$ stands for $g(b)-g(a)$.
Denote $a(s)=\int_{\rho}^{s}g_{n}(x-t)\, dt$, and remember that (\ref{eq:intgparts})
gives that $|a|\leq C$. Hence when we  integrate by parts the integral on
the right hand side of (\ref{eq:Int01notI}) gives \begin{align*}
\left|\int_{\mathbb{T}\setminus I}H^{(D)}(t)g_{n}(x-t)\, dt\right| & \leq\left|H^{(D)}(-\rho)a(2\pi-\rho)\right|+\left|\int_{\rho}^{2\pi-\rho}H^{(D+1)}(t)a(t)\, dt\right|\\
 & \stackrel{(\ref{eq:intgparts})}{\leq}C\left|H^{(D)}(-\rho)\right|+C\int_{\rho}^{2\pi-\rho}\left|H^{(D+1)}(t)\right|dt\\
 & \stackrel{(\ref{eq:Hilbert})}{\le}\frac{(CD)^{CD}}{\rho^{D+1}}.\end{align*}
Similarly we can use (\ref{eq:fstarderivd}) and (\ref{eq:Hilbert})
to estimate the sum in (\ref{eq:Int01notI}) and this gives
\begin{equation}
\left|\left(\frac{d}{dx}\right)^{D}\int_{\mathbb{T}\setminus I}H(t)g_{n}(x-t)\, dt\right|\leq\frac{(CD)^{CD}}{\rho^{D+1}}.\label{eq:int01notIf}\end{equation}
Together with (\ref{eq:intIf}) and (\ref{eq:intIint1mI}) this gives
\begin{equation}
\left|\widetilde{g_{n}}^{(D)}(x)\right|\leq\frac{(CD)^{CD}}{\rho^{D+1}}\quad\forall x\not\in K_{n}\cup K'.\label{eq:gntildeD}\end{equation}
 Any $x\not\in K'$ is also not in some $K_{m}$ and hence (\ref{eq:gntildeD})
holds for any $n>m$ and hence it holds for $\widetilde{g}$. With
(\ref{eq:fstarderivd}) and (\ref{eq:forlem1}) we can use Lemma \ref{lem:fsmthFsmth}
and get (\ref{eq:Fstarsmooth}) which proves Lemma \ref{lem:Fstarsmooth}.
\hfill\qed

\Subsec{\label{sub:Taylor}The Taylor coefficients of $F$}
 In this section we prove Lemma \ref{lem:Taylor0}, namely
show that with probability 1, $\widehat{F}(m)\rightarrow0$ as $m\rightarrow\infty$.
First we define $F_{n}$ to be the harmonic extension of $f_{n}=e^{g_{n}+i\widetilde{g_{n}}}$
to $\mathbb{D}$ (each $f_{n}$ is bounded) and we want to find some
$n$ such that $\widehat{f_{n}}(m)=\widehat{F_{n}}(m)$ approximates
$\widehat{F}(m)$. Summing (\ref{eq:fnDfn1DdzK}) and (\ref{eq:zntilde})
over $n$ we get\[
|(G_{n}+i\widetilde{G_{n}})(z)-(G+i\widetilde{G})(z)|\leq\frac{C}{(1-|z|)2^{n}}.\]
Fix, therefore, $n=n(m):=\left\lceil C\log m\right\rceil $ for some
$C$ sufficiently large, and get, for every $z$ with $|z|=1-\frac{1}{m}$
that $|(G_{n}+i\widetilde{G_{n}})(z)-(G+i\widetilde{G})(z)|\leq1/m$.
Further, we have that \begin{equation}
\max\left|F_{n}\right|\stackrel{(\ref{eq:maxfn})}{\leq}e^{o(n)}\leq m^{o(1)}\label{eq:maxFnmo1}\end{equation}
which means that, for $|z|=1-\frac{1}{m}$, $$
|F_{n}(z)-F(z)|\leq|F_{n}(z)\Vert 1-\exp((G_{n}+i\widetilde{G_{n}})(z)-
(G+i\widetilde{G})(z))|\leq Cm^{-1/2}.$$
Finally we use $$
\widehat{F}(m)=\frac{1}{m!}F^{(m)}(0)=\int_{|z|=1-1/m}z^{-m-1}F(z)\, dz
\pagebreak $$
so that
\begin{equation}
|\widehat{F_{n}}(m)-\widehat{F}(m)|=\left|\int_{|z|=1-1/m}z^{-m-1}(F_{n}(z)-F(z))\, dz\right|\leq Cm^{-1/2}\label{eq:FnFmhalf}\end{equation}
and we see that it is enough to calculate $\widehat{f_{n}}(m)$.

\Subsec{\label{sub:Probability}Probability}
At this point   we use the fact that the $s(n,k)$ are random.
Take them to be independent and uniformly distributed on $[0,1]$.
We shall perform a (probabilistic) estimate of $\widehat{f_{n}}(m)$
by moment methods. Unfortunately, it seems we need the fourth moment.
We start with a lemma that contains the calculation we need without
referring to analytic functions

For $i=1,2,3,4$ and $j=1,2,3$ we denote $i\subset j$ if $i=j$
or $i=4$ and $j=3$. The inverse will be denoted by $i\not\subset j$.

\begin{lem}
\label{lem:real}Let $I_{i}$ be $4$ intervals and let $\tau,\alpha,\beta>0$
be some numbers. Let $h_{1},h_{2},h_{3}$ be functions satisfying
\begin{align}
\int_{I_{i}}|h_{j}| & \leq\alpha, & i & \subset j,\label{eq:estGiIi}\\
|h_{j}(x)| & =1,|h_{j}'(x)|\leq\beta,|h_{j}''(x)|\leq\beta^{2}\quad\forall x\in I_{i}+[-\tau,\tau] ,& i &
\not\subset j,\label{eq:estGtg}\end{align} 
where \/{\rm ``}\/$+$\/{\rm ''}\/ stands for regular set addition. Let $t_{1}$ and
$t_{2}$ be two random variables{\rm ,} uniformly distributed on $[0,\tau]${\rm ,}
and let $t_{3}=t_{4}=0$. Define\begin{equation}
f(x)=f_{t_{1},t_{2}}(x):=h_{1}(x-t_{1})h_{2}(x-t_{2})h_{3}(x).\label{eq:defF}\end{equation}
Then\begin{equation}
E:=\left|\mathbb{E}\left(\prod_{i=1}^{4}\int_{I_{i}+t_{i}}f(x_{i})e^{-imx_{i}}\, dx_{i}\right)\right|\leq C\frac{\alpha^{4}}{m^{2}}\left(\max\beta,\frac{1}{\tau}\right)^{2}.\label{eq:defE}\end{equation}
\end{lem}
\vskip8pt

\Proof 
Denote $\beta'=\max\beta,\frac{1}{\tau}$. Define $S_{i}=\int_{I_{i}+t_{i}}f(x_{i})e^{-imx_{i}}\, dx_{i}$.
Translate $S_{1}$ and $S_{2}$ to get\begin{alignat*}{1}
S_{1} & =\int_{I_{1}}h_{1}(x_{1})h_{2}(x_{1}+t_{1}-t_{2})h_{3}(x_{1}+t_{1})
e^{-im(x_{1}+t_{1})}\, dx_{1},\\
S_{2} & =\int_{I_{2}}h_{1}(x_{2}+t_{2}-t_{1})h_{2}(x_{2})h_{3}(x_{2}+t_{2})e^{-im(x_{2}+t_{2})}\, dx_{2}.\end{alignat*}
Changing the order of integration we get\begin{equation}
E=\left|\int_{I_{1}}h_{1}(x_{1})e^{-imx_{1}}\dotsi\int_{I_{3}}h_{3}(x_{3})e^{-imx_{3}}\int_{I_{4}}h_{3}(x_{4})e^{-imx_{4}}A\, dx_{1}\dotsm dx_{4}\right|\label{eq:EU1U4A}\end{equation}
where $A$, the central element, is defined by\[
A=\frac{1}{\tau^{2}}\int_{0}^{\tau}\int_{0}^{\tau}e^{-im(t_{1}+t_{2})}A(t_{1},t_{2})\, dt_{1}\, dt_{2}\]
and where \begin{align*}
A(t_{1},t_{2}):=\; & h_{2}(x_{1}+t_{1}-t_{2})h_{3}(x_{1}+t_{1})h_{1}(x_{2}+t_{2}-t_{1})h_{3}(x_{2}+t_{2})\\
 & h_{1}(x_{3}-t_{1})h_{2}(x_{3}-t_{2})h_{1}(x_{4}-t_{1})h_{2}(x_{4}-t_{2}).\end{align*}
 The lemma will be proved once, we estimate $A$, which will be done
by integrating by parts over $t_{1}$ and then over $t_{2}$. We notice
that $A$ contains only expressions of the type $h_{j}(x)$ where
$x\in I_{i}+[-\tau,\tau]$ and $i\not\subset j$. Therefore, using
(\ref{eq:estGtg}) we get\[
|A(t_{1},t_{2})|=1,\quad\left|\frac{\partial A(t_{1},t_{2})}{\partial t_{i}}\right|\leq5\beta,\quad\left|\frac{\partial^{2}A(t_{1},t_{2})}{\partial t_{1}\partial t_{2}}\right|\leq25\beta^{2}.\]
 Integrating by parts once, we get\[
\left|\int_{0}^{\tau}A(t_{1},t_{2})e^{-imt_{1}}\, dt_{1}\right|\leq\frac{2}{|m|}+\frac{5\tau\beta}{|m|}\leq\frac{7\tau\beta'}{|m|}.\]
Further,
\begin{eqnarray*}
\left|\frac{\partial}{\partial t_{2}}\int_{0}^{\tau}
A(t_{1},t_{2})e^{-imt_{1}}\, dt_{1}\right| & 
=&\left|\int_{0}^{\tau}\frac{\partial A(t_{1},t_{2})}{\partial
t_{2}}e^{-imt_{1}}\, dt_{1}\right|\\
 & =&\Bigg|\frac{\partial A(t_{1},t_{2})}{\partial t_{2}}
\frac{e^{-imt_1}}{-im}\Big|_{0}^{\tau}-\int_{0}^{\tau}\frac{\partial^{2}A(t_{1},t_{2})}{\partial
t_{2}\partial t_{1}}\frac{e^{-imt_1}}{-im}\Bigg|\\
&  \leq & \frac{10\beta}{|m|}+\frac{25\tau\beta^{2}}{|m|}\leq\frac{35\tau\beta^{'2}}{|m|}.
\end{eqnarray*}
These two statements  allow us to perform integration by parts over $t_{2}$ getting\[
|A|=\frac{1}{\tau^{2}}\left|\int_{0}^{\tau}\int_{0}^{\tau}A(t_{1},t_{2})e^{-im(t_{1}+t_{2})}\right|\leq\frac{1}{\tau^{2}}\left(\frac{14\tau
\beta'}{m^{2}}+\frac{35\tau^{2}\beta^{'2}}{m^{2}}\right)\leq\frac{49\beta^{'2}}{m^{2}}.\]
Plugging this into (\ref{eq:EU1U4A}) and integrating using (\ref{eq:estGiIi}) we
conclude the proof of  the lemma.
\Endproof\vskip4pt 
 
Continuing the proof of the theorem, for every $0\leq k<2^{n}$ denote
\[
\mathcal{I}_{k}=\int_{I(n,k)}f_{n}(x)e^{-imx}\, dx.\]
We note that (\ref{eq:maxfn}) shows that \begin{equation}
|\mathcal{I}_{k}|\leq\int_{I(n,k)}|f_{n}(x)|\leq\sigma_{n}e^{o(n)}=:\gamma.\label{eq:defalpha}\end{equation}
In other words, $\gamma=\gamma(n)$ is a bound for $|\mathcal{I}_{k}|$
independent of $k$ satisfying \begin{equation}
\gamma=\sigma_{n}e^{o(n)}=2^{-n}m^{o(1)}\label{eq:gammanmo1}\end{equation}
(for the last equality, remember that $\sigma_{n}=2^{-n}\Phi(n)$,
$\Phi(n)=n^{-o(1)}$ (\ref{eq:Phinsmall}) and $n\approx\log m$).

\begin{lem}
\label{lem:complex}Let $0\leq k_{1},k_{2},k_{3},k_{4}<2^{n}$ and
let $1\leq r<n${\rm ,} and assume that the $I(n,k_{i})$ belong to at least
three different intervals of rank $r$. Then \[
\mathbb{E}(\mathcal{I}_{k_{1}}\mathcal{I}_{k_{2}}\mathcal{I}_{k_{3}}\mathcal{I}_{k_{4}})\leq\gamma^{4}\frac{C\omega(n)^{2}}{m^{2}\tau_{r}^{3}}.\]
\end{lem}

\Proof 
Define $q_{1},\dotsc,q_{4}$ using $I(n,k_{i})\subset I(r,q_{i})$.
We may assume without loss of generality that the two $q_{i}$-s which
may be equal are $q_{3}$ and $q_{4}$. Let $\mathcal{X}$ be the
$\sigma$-field spanning all $s$-es \emph{except} $s(r,q_{1})$ and
$s(r,q_{2})$. We shall show\[
\mathbb{E}(\mathcal{I}_{k_{1}}\mathcal{I}_{k_{2}}\mathcal{I}_{k_{3}}\mathcal{I}_{k_{4}}|\mathcal{X})\leq\gamma^{4}\frac{C\omega(n)^{2}}{m^{2}\tau_{r}^{3}}\]
and then integrating over $\mathcal{X}$ will give the result. We
note that conditioning by $\mathcal{X}$ is in effect fixing everything
except the positions of $I(r,q_{1})$ and $I(r,q_{2})$ inside $I(r-1,\left\lfloor q_{i}/2\right\rfloor )$.
To be more precise, two copies of $l$ also move with $I(r,q_{i})$.
Therefore define $J_{j}:=I(r,q_{j})+\left[-\tau_{r},\tau_{r}\right]$
($j=1,2$), which is the part of $f_{n}$ that moves when $s(r,q_{j})$
changes (there are zones where $f_{n}\equiv-\mu(r)$ which expand
and contract on the sides of $J_{j}$) and denote $J_{3}=\mathbb{T}\setminus(J_{1}\cup J_{2})$.
Assume for a moment that $s(r,q_{1})=s(r,q_{2})=0$ and define, using
this assumption,\begin{equation}
\begin{aligned}\eta_{j} & :=(g_{n}+\omega(r))|_{J_{j}},\quad  & j & =1,2,3, & h_{j} &
:=e^{\eta_{j}+i\widetilde{\eta_{j}}},\\ I_{i} & :=I(n,k_{i}),\quad  & i & =1,2,3,4.\end{aligned}
\label{eq:defIi}\end{equation}
Under the assumption $s(r,q_{1})=s(r,q_{2})=0$ we clearly have $f_{n}=h_{1}h_{2}h_{3}e^{-\omega(r)}$
and when we remove this assumption, the only change is a translation
of $h_{1}$ and $h_{2}$. In other words, if we define $t_{i}=s(r,q_{i})\tau_{r}$
then $$f_{n}(x)=h_{1}(x-t_{1})h_{2}(x-t_{2})h_{3}(x)e^{-\omega(r)}.$$
Examining (\ref{eq:defF}) we see that $|\mathbb{E}(\mathcal{I}_{k_{1}}\mathcal{I}_{k_{2}}\mathcal{I}_{k_{3}}\mathcal{I}_{k_{4}}|\mathcal{X})|=e^{-4\omega(r)}E$
where $E$ is defined by (\ref{eq:defE}); where the $I_{i}$ of (\ref{eq:defE})
are the same as those of (\ref{eq:defIi}); and where the $\tau$
of (\ref{eq:defE}) is $\tau_{r}$. To make (\ref{eq:defE}) concrete
we need to specify values for the $\alpha$ and $\beta$ of (\ref{eq:estGiIi})
and (\ref{eq:estGtg}) and prove that they hold. We define\[
\alpha=\gamma e^{\omega(r)},\quad\beta=C\frac{\omega(n)}{\tau_{r}^{3/2}}.\]
Notice that $\beta$ is obviously larger than $1/\tau_{r}$. With
all these, Lemma \ref{lem:complex} would follow from Lemma \ref{lem:real}
once we show (\ref{eq:estGiIi}) and (\ref{eq:estGtg}). (\ref{eq:estGiIi})
is clear from the definitions of $\alpha$ above, $\gamma$ (\ref{eq:defalpha})
and $h_{i}$ (\ref{eq:defIi}), so we need only show (\ref{eq:estGtg}).

Examining the definitions of $\eta_{j}$ and $I(n,k)$ we see easily that $\eta_{j}(x)=0$ for $x\in
I_{i}+[-2\tau_{r},2\tau_{r}]$ whenever $i\not\subset j$ (we defined $l^{\pm}$ (\ref{eq:deflpm})
with a slightly larger ``space'' so that this fact would be true).
This immediately shows $|h_{i}(x)|=1$. Further, $h_{i}'=h_{i}(\eta_{i}'+i\widetilde{\eta_{i}}')$
gives $|h_{i}'|=|\widetilde{\eta_{i}}'|$ and $h_{i}''=h_{i}((\eta_{i}'+i\widetilde{\eta}_{i}')^{2}+\eta_{i}''+i\widetilde{\eta_{i}}'')$
gives $|h_{i}''|\leq|\widetilde{\eta_{i}}'|^{2}+|\widetilde{\eta_{i}}''|$.
As in (\ref{eq:Int01notI}), the derivatives of $\widetilde{\eta_{i}}$
have the representations \[
\widetilde{\eta_{i}}'(x)=\int_{\mathbb{T}}\eta_{i}(x-t)H'(t)\,
dt, \quad\widetilde{\eta_{i}}''(x)=\int_{\mathbb{T}}\eta_{i}(x-t)H''(t)\, dt\] where $H$ is the Hilbert kernel.
In general there are boundary terms (as in the calculation in (\ref{eq:Int01notI})), but, as remarked,
in our case $\eta_{i}$ is zero in $[x-\tau_{r},x+\tau_{r}]$ (when
$x\in I_{j}+[-\tau_{r},\tau_{r}]$) so these terms disappear. We may
therefore estimate\begin{equation}
\widetilde{\eta_{i}}'(x)\leq\Vert \eta\Vert _{2}\left\Vert H'|_{\mathbb{T}\setminus[-\tau_{r},\tau_{r}]}
\right\Vert _{2},\qquad\widetilde{\eta_{i}}''(x)\leq\Vert \eta\Vert _{2}\left\Vert
H''|_{\mathbb{T}\setminus[-\tau_{r},\tau_{r}]}\right\Vert _{2}.\label{eq:utilde}\end{equation} The
second terms are a straightforward calculation from (\ref{eq:Hilbert}) and we get \[
\left\Vert H'|_{[-\tau_{r},\tau_{r}]^{c}}\right\Vert _{2}\approx\tau_{r}^{-3/2},\qquad\left\Vert
H''|_{[-\tau_{r},\tau_{r}]^{c}}\right\Vert _{2}\approx\tau_{r}^{-5/2}.\] The first terms can be estimated
easily: since the singularities in
$\eta$ (which originally came from the $x^{-1/3}$ factor in $l$)
are in $L^{2}$. Indeed, it is for this point that we defined $l$
using $x^{-1/3}$. We easily get \begin{equation}
\Vert \eta\Vert _{2}\leq C\omega(n)\label{eq:normg}\end{equation}
which gives us the estimate we need:\[
|h_{i}'|\leq C\tau_{r}^{-3/2}\omega(n),\qquad|h_{i}''|\leq C\tau_{r}^{-3}\omega^{2}(n).\]
With this the conditions of Lemma \ref{lem:real} are fulfilled and
we are done.
\hfill\qed

\demo{Proof of Lemma {\rm \ref{lem:Taylor0}}} Define \[
X=X_{m}=\sum_{k=0}^{2^{n}-1}\int_{I(n,k)}f_{n}(x)e^{-imx}\, dx.\]
Now the difference $\widehat{f_{n}}(m)-X$ is exactly $\widehat{f_{n}\mathbf{1}_{\mathbb{T}\setminus K_{n}}}(m)$.
The functions $f_{n}\mathbf{1}_{\mathbb{T}\setminus K_{n}}$ are uniformly
$C^{1}$ so this difference converges to zero. To see this last claim,
note that (\ref{eq:fstarderivd}) and (\ref{eq:gntildeD}) show that
$g_{n}'+i\left(\widetilde{g_{n}}\right)'$ has a bound of $C/d(x,K_{n})^{2}$
while (\ref{eq:forlem1a}) and  (\ref{eq:forlem1}) show that $f_{n}\leq Cd(x,K_{n})^{10}$
(actually $f_{n}$ converges to zero near $K_{n}$ superpolynomially
uniformly).

Therefore we want to bound $X$, and we shall estimate $\mathbb{E}X^{4}$.
Let\[
E(k_{1},k_{2},k_{3},k_{4}):=\mathbb{E}\prod\mathcal{I}_{k_{i}}; \]
 let $r(k_{1},\dotsc,k_{4})$ be the minimal $r$ such that the $I(n,k_{i})$-s
are contained in at least $3$ distinct intervals of rank $r$. A
simple calculation shows \[
\#\{(k_{1},\dotsc,k_{4}):r(k_{1},\dotsc,k_{4})=r\}\approx2^{4n-2r}.\]

If $\tau_{r}$ is too small then the estimate of the lemma is useless
and it would be better to estimate $|E(k_{1},\dotsc,k_{4})|\leq\gamma^{4}$.
Let $R$ be some number. Then for large $r$ we have the estimate\begin{equation}
E_{1}:=\sum_{r(k_{1},\dotsc,k_{4})\geq R}E(k_{1},\dotsc,k_{4})\leq C\gamma^{4}2^{4n-2R}\stackrel{(\ref{eq:gammanmo1})}{=}m^{o(1)}2^{-2R}.\label{eq:E1}\end{equation}

For small $r$ we use the lemma to get a better estimate. Examine
one such $k_{1},\dotsc,k_{4}$ and let $r=r(k_{1},\dotsc,k_{4})$.
Lemma \ref{lem:complex} gives\[
E(k_{1},\dotsc,k_{4})\leq\gamma^{4}\frac{C\omega^{2}(n)}{m^{2}\tau_{r}^{3}}\stackrel{(*)}{=}\frac{\gamma^{4}}{m^{2-o(1)}\tau_{r}^{3}}\]
where $(*)$ comes from the regularity condition $\omega(n)=e^{o(n)}$
and $n\approx\log m$. Therefore
\begin{eqnarray}
E_{2}&:  =&\sum_{r(k_{1},\dotsc,k_{4})<R}E(k_{1},\dotsc,k_{4})\leq
\gamma^{4}2^{4n}m^{-2+o(1)}\sum_{r=1}^{R}2^{-2r}\tau_{r}^{-3}
\label{eq:E2}\\
&\stackrel{(\ref{eq:taunreg},\ref{eq:gammanmo1})}{=}&
m^{-2+o(1)}\sum_{r=1}^{R}2^{r+o(r)}=m^{-2+o(1)}2^{R+o(R)}.\nonumber
\end{eqnarray} Taking $R=\left\lfloor \frac{2}{3}\log_{2}m\right\rfloor $ and summing (\ref{eq:E1}) and
(\ref{eq:E2}) we get\begin{equation}
\mathbb{E}X^{4}\leq m^{-4/3+o(1)}.\label{eq:EX4m43}\end{equation}
This gives that \[
\mathbb{E}\sum_{m}X_{m}^{4}\leq\sum\mathbb{E}X_{m}^{4}<\infty.\]
In particular, with probability $1$, $X_{m}\to0$. As remarked above,
this shows that $\widehat{f_{n}}(m)\rightarrow0$ and hence $\widehat{F}(m)\to0$
which concludes Lemma \ref{lem:Taylor0} and the theorem. \phantom{overthere}
\hfill\qed

\Subsec{Remarks}
 1. It is clear that if $f\in\PLA$ then the associated analytic
function $F$ defined by (\ref{eq:zPLA}), (\ref{eq:defPLA}) has
the estimate\begin{equation}
F(z)=o\left(\frac{1}{1-|z|}\right)\label{eq:Fess}\end{equation}
simply because $\widehat{F}(n)\to0$. It turns out that in some vague
sense, this inequality is the ``calculationary essence'' of $\PLA\setminus H^{2}$.
In other words, if you have a singular distribution whose analytic
part $F$ satisfies (\ref{eq:Fess}) and its boundary value is in
$L^{2}$ then you are already quite close to constructing a nonclassic
$\PLA$ function. Note that (\ref{eq:Fess}) is enough to prove uniqueness
(see Theorem \ref{thm:unique}$'$ below) and the additional information
$\widehat{F}(n)\to0$ does not help. This ideology also stands behind
the proof above. To understand why, let $K$ be a nonprobabilistic
Cantor set with the same thickness; namely at step $n$ the total
length of the $2^{n}$ intervals is $\Phi(n)$. Let $G$ be a harmonic
function constructed similarly; i.e.\ ``hang'' copies of $-x^{-1/3}$
suitably dilated and shifted from all intervals contiguous to $K$
and define $F=e^{G+i\widetilde{G}}$. Then a much simpler calculation
shows that $F$ satisfies (\ref{eq:Fess}), and even the stronger\begin{equation}
F(z)=\left(\frac{1}{1-|z|}\right)^{o(1)}.\label{eq:essence2}\end{equation}
 This implies\[
\sum_{k=2^{n}}^{2^{n+1}}|\widehat{F}(k)|^{2}=o(2^{n})\]
so that ``in average'' the coefficients tend to zero, with no need
for probability in the construction. Thus the probabilistic skewing
introduced above ``smears'' the spectrum of $F$ uniformly and
allows us  to conclude the stronger $\widehat{F}(m)\rightarrow0$. 

\demo{Question}
Is the nonprobabilistic construction (e.g.~taking all $s(n,k)$
to be $0$) in $\PLA$?
\Enddemo

The use of stochastic perturbations of the time domain to smooth
singularities in the spectrum is not new. One may find examples of
such techniques in \cite{K85}, notably the use of Brownian images
in Chapter 17, and in \cite{KO98}. A reader fluent in these techniques
would probably assume it is possible to simplify the proof of Lemma
\ref{lem:Taylor0} significantly along the following rough lines:
find some event $\mathcal{X}$ that would separate $\mathcal{I}_{k_{1}}$
from the rest of the $\mathcal{I}_{k_{i}}$'s making them independent, 
perhaps similar to the $\mathcal{X}$ actually used. Now calculate
$\mathbb{E}(\mathcal{I}_{k}|\mathcal{X})$ using one simple integration
by parts. Next multiply $\mathbb{E}(\mathcal{I}_{k_{i}}|\mathcal{X})$
out, and integrate over $\mathcal{X}$. Unfortunately, it seems that
no proof can be constructed this way. The problem is that, while $g$
has a local structure and would be amenable to such a handling, $\widetilde{g}$
does not, and any change to one $s(n,k)$ affects the values of $\widetilde{g}$
globally. 

\smallskip{}
  2. The regularity condition $\omega(n)=e^{o(\log n)}$ can
be relaxed somewhat, but it is not clear whether it can be removed
altogether. For example, there is an inherent difficulty in generalizing
Lemma \ref{lem:fsmthFsmth} without this condition.

\smallskip{}
  3. It is also of interest to ask how fast does $\widehat{F}(m)\to0$,
or in other words how much do we pay in $\{ c(n)\}_{n>0}$ for the
quick decrease of the $\{ c(n)\}_{n<0}$. Using Chebyshev's inequality
with (\ref{eq:EX4m43}) it is easy to see that $|\widehat{F}(m)|\leq m^{-1/12+o(1)}$.
This, however, can be improved significantly. Indeed, one may change
the definition of $l$, (\ref{eq:lsing}) to have a singularity of
type $x^{-\varepsilon}$ and then replace (\ref{eq:utilde}) with an
$L^{p}-L^{q}$ estimate and get $\tau_{r}^{2+\varepsilon}$ in the formulation
of Lemma \ref{lem:complex}, which would end up as $X_{m}\leq m^{-1/4+\varepsilon}$
almost surely. Further, it is possible to use higher moment  estimates.
To estimate the $2k^{\rm th}$ moment, use a generalized version of Lemmas
\ref{lem:real} and \ref{lem:complex} for $k$ moving intervals and
$k$ stationary ones to get an estimate of $m^{-k}\tau_{r}^{k+\varepsilon}$
and the final outcome would be $X_{m}\leq m^{-1/2+1/2k+\varepsilon}$
almost surely. Thus in effect we may construct a function $F$ satisfying
$\widehat{F}\in l^{2+\varepsilon}$ for all $\varepsilon>0$, almost surely. 

\smallskip{}
  4. It is possible to characterize precisely the size of
exceptional sets for the ``nonclassic'' part of $\PLA\cap\, L^{2}$.
Namely, denote by $\Lambda$ the (generalized) Hausdorff measure generated
by the function $t\mapsto t\log1/t$. Then the following is true

\medbreak \hangindent=29pt\hangafter=1 {\scshape Theorem.} 
 {\rm (i)}  {\it There exists a function $f\in L^{2}\setminus H^{2}$
which admits a decomposition {\rm (\ref{eq:defPLA})} converging everywhere
outside of some compact $K$ of finite $\Lambda$-measure.}

\vskip4pt \noindent \hglue4pt
\hangindent=28pt\hangafter=1 {\rm  (ii)}\hskip9pt {\it The result fails if one replaces the condition
$\Lambda(K)<\infty$ with $\Lambda(K)=0$ even if $K$ is not required
to be compact.}
 \vskip4pt

Part (i) can be  proved by the construction of the section,
with a nonnegligible simplification since we do not watch for smoothness.
Part (ii) follows easily from Phragm\'en-Lindel\"of-like
theorems for analytic functions of slow growth in~$\mathbb{D}$. See
\cite{B92}, \cite{D77}.

Note that in the symmetric settings, the corresponding exceptional
sets (so call\-ed $M$-sets) could have dimension zero \cite{B64}, \cite{KS94},
\cite{KL87}
and moreover, may be ``thin'' with respect to any (generalized)
Hausdorff measure \cite{I68}.

\smallskip{}
 5. We also have some structural results about the set $\PLA\cap\, C$.
Namely:

\vskip4pt {\scshape Theorem}.
{\it $\PLA\cap\, C$ is the  first Baire category and has zero Wiener measure.}

\vskip4pt {\scshape Theorem}.
{\it $\PLA$ is dense in $C$ \/{\rm (}\/in sharp contrast to $H^{2}${\rm )}.}
\vskip6pt

We intend to publish proofs of these three results elsewhere.\footnote{The last two results are to appear at
{\it Bull.\ London Math.\ Soc.}}

\vglue-20pt
\phantom{up}
\section{Uniqueness}
\vglue-4pt

4.1.\quad The most natural settings for the statement of the uniqueness result
is that of boundary behavior of analytic functions. Let us therefore
restate Theorem \ref{thm:unique} in a stronger form 

\demo{\scshape Theorem \ref{thm:unique}$'$} {\it  Let $F$ be an analytic function on $\mathbb{D}$
satisfying\begin{equation} F(z)=O\left(\frac{1}{(1-|z|)^{M}}\right)\quad\textrm{for some
}M.\label{eq:Fpoly}\end{equation} Assume that $F$ has nontangential boundary limit almost everywhere
and that \begin{equation}
F(\zeta)=f(\zeta):=\sum_{n=-\infty}^{-1}c(n)\zeta^{n}\quad\textrm{a.e. on }\partial\mathbb{D}\label{eq:defcnneg}\end{equation}
and assume the $c(n)$ satisfy {\rm (\ref{eq:qa})}
 with some $\omega:\mathbb{R}^{+}\to\mathbb{R}^{+}${\rm ,}
$\omega(t)/t$ increasing and $\sum\frac{1}{\omega(n)}<\infty$. Then
$F$ and $f$ are identically zero.}
\Enddemo

To see that Theorem \ref{thm:unique}$'$ generalizes Theorem \ref{thm:unique}
define $F(z)$ by (\ref{eq:zPLA}) and note that (\ref{eq:Fess})
is stronger than (\ref{eq:Fpoly}). And as usual, Abel's theorem shows
that $F$ has nontangential boundary limit a.e.\

In this section the notation $C$ and $c$ will be allowed to depend
on the function $F$ --- here we consider $F$ as given and fixed.
By $\mathbf{m}$ we denote the arc length on the circle, normalized
so that $\mathbf{m}\partial\mathbb{D}=1$. For $\theta\in\partial\mathbb{D}$
we denote by $I(\theta,\varepsilon)$ an arc centered around $\theta$
with $\mathbf{m}I(\theta,\varepsilon)=\varepsilon$.

For a compact subset $E$ of $\partial\mathbb{D}$ we shall define
the \emph{Privalov domain} over $E$, $\mathcal{P}(E)$, to be a subset
of $\mathbb{D}$ created by removing, for every arc $I$ from the
complement of $E$, a disk $D_{I}$ orthogonal to $\mathbb{D}$ at
the end points of $I$. If $I$ is larger than a half circle, remove
$\mathbb{D}\setminus D_{I}$ instead of $D_{I}$ so that in both cases
you remove the component containing $I$ (this definition is slightly
different from the standard one). The following is well known:

\begin{lem}
\label{lem:privdom} Let $F$ be an analytic function on $\mathbb{D}$
with almost everywhere nontangential boundary values{\rm ,} and let $\delta>0$.
Then there exists a compact set $E\subset\partial\mathbb{D}$ with
$\mathbf{m}E>1-\delta$ such that $F$ is continuous on $\mathcal{P}(E)$.
\end{lem}

\demo{{\rm 4.2}} The following lemma is simple but plays a crucial part in the proof. 

\begin{lem}
\label{lem:main}Let $L$ be a function on some measure space with
a probability measure $\mu$ with $A\leq L\leq B$. Assume for some
$\varepsilon\in[0,1]${\rm ,}
\begin{equation}
\int L\, d\mu=\varepsilon B+(1-\varepsilon)A.\label{eq:intL}\end{equation}
Let $D\in\left]A,B\right[$. Then\begin{equation}
\int\max\{ L,D\}\, d\mu\leq\varepsilon B+(1-\varepsilon)D.\label{eq:maxLD}\end{equation}
\end{lem}

\Proof 
By considering $L-D$, we may assume without loss of generality that
$D=0$. Assume by contradiction that \[
\int L^{+}>\varepsilon B.\]
This shows that the support of the positive part of $L$ has measure
$>\varepsilon$ so that the support of the negative part must have measure
$<1-\varepsilon$ and therefore $\int L^{-}>A(1-\varepsilon)$, a contradiction
to (\ref{eq:intL}).
\hfill\qed

\begin{lem}
\label{lem:harmpriv}Let $E$ be   compact in $\partial\mathbb{D}$
and let $\varepsilon>0$. Let $z\in\mathcal{P}(E)$ with $|z|>1-\varepsilon$
and define $\zeta:=z/|z|$ and $l:=(1-|z|)/\varepsilon$. Assume that\begin{equation}
\mathbf{m}(I(\zeta,l)\setminus E)\leq\varepsilon^{2}l.\label{eq:Izlemp}\end{equation}
Let $\mathcal{D}$ be the component of $\mathcal{P}(E)\setminus(1-l)\mathbb{D}$
containing $z$. Then the harmonic measure of $E$ has the estimate\begin{equation}
\Omega(z,\mathcal{D})(E)>1-C_{1}\varepsilon.\label{eq:lemharmpriv}\end{equation}
\end{lem}
 \vskip8pt

Here $C_{1}$ is an absolute constant. Similarly all constants in
the proof of the lemma, both explicit and implicit in $\approx$ notations
are absolute.

\Proof 
Without loss of generality, assume $\varepsilon<\frac{1}{4}$. For any
arc $I$ which is a component of $\partial\mathbb{D}\setminus E$
let $D_{I}$ be the Privalov disk as in the definition of a Privalov
domain. Let $B$ be a Brownian motion starting from $z$. Let $T$
be the stopping time of $B$ on $\partial\mathbb{D}$ and   $T_{I}$
be the stopping time of $B$ on $D_{I}$.

We start with an estimate of $p:=\Omega(z,\mathbb{D}\setminus D_{I})(\partial D_{I})$.
We need to consider two cases: 
\smallbreak
(i)  If $I$ is ``far'' from $z$ we use $p\leq C(1-|z|)/d(z,I)$.

\smallbreak (ii) \label{enu:Iclose}If $I$ is ``close'' to $z$ we use $p\leq C\mathbf{m}I/(1-|z|)$.
\smallbreak\noindent 
Both follow from the conformal invariance of the harmonic measure
\cite[Theorem V.1.2]{B95} which gives an explicit formula for $\Omega(z,\mathbb{D}\setminus D_{I})$
and (i) and (ii) with a simple calculation.

Let $J$ be $\mathbb{T}\setminus I(\zeta,l/2)$, and let $T_{J}$
be the hitting time of $J$. Then (i) shows that $\Omega(z,\mathbb{D}\setminus D_{J})(\partial D_{J})\leq C\varepsilon$
and hence\[
\mathbb{P}(T_{J}<T)\leq C\varepsilon.\]
For any $I\subset J$ we have $D_{I}\subset D_{J}$ so that  $T_{I}>T_{J}$.
Now,  \begin{equation}
\mathbb{P}\left(\inf_{I\subset J}T_{I}<T\right)\leq C\varepsilon.\label{eq:IinJ}\end{equation}

Next, if $I\subset I(\zeta,l)$ we use (ii) above and
(\ref{eq:Izlemp}) to get \begin{equation}
\mathbb{P}\left(\inf_{I\subset I(\zeta,l)}T_{I}<T\right)\leq\sum_{I\subset I(\zeta,l)}\mathbb{P}(T_{I}<T)\stackrel{\textrm{(ii)}}{\leq}C\sum_{I\subset I(\zeta,l)}\frac{\mathbf{m}I}{1-|z|}\stackrel{(\ref{eq:Izlemp})}{\leq}C\varepsilon.\label{eq:IninJ}\end{equation}
The assumption $\varepsilon\leq\frac{1}{4}$ together with (\ref{eq:Izlemp})
shows that any arc $I$ in the complement of $E$ is either a subset
of $I(\zeta,l)$ or of $J$. Hence the lemma is almost finished. We
still have to deal with $\{|z|=1-l\}$, which is easy: define $T^{*}$
to be the stopping time of $B$ on the circle $\{|z|=1-l\}$. Let
$\mathcal{A}$ be the annulus $\{1-l\leq|z|\leq1\}$. The harmonic
measure $\Omega(z,\mathcal{A})(\{|z|=1-l\})$ can be calculated explicitly
from the fact that the solution $h$ of Dirichlet's problem on $\mathcal{A}$
with boundary conditions $1$ on $\{|z|=1-l\}$ and $0$ on $\{|z|=1\}$
has the explicit form\begin{equation}
h(w)=\frac{\log|w|}{\log(1-l)}.\label{eq:harmann}\end{equation}
Since $\Omega(z,\mathcal{A})(\{|z|=1-l\})\stackrel{(\ref{eq:kak})}{=}h(z)\leq C\varepsilon$,
 \[
\mathbb{P}(T^{*}<T)\leq C\varepsilon\]
and this  together with (\ref{eq:IinJ}) and (\ref{eq:IninJ}) gives
\[
\mathbb{P}\left(\min\left\{ \inf_{I}T_{I},T^{*}\right\} <T\right)\leq C\varepsilon.\]
This is equivalent to (\ref{eq:lemharmpriv}) and the lemma is proved.
\hfill\qed

\setcounter{Subsec}{2}

\Subsec{Proof of Theorem {\rm \ref{thm:unique}}$'$}
Without loss of generality we may assume the $\log$ in (\ref{eq:qa})
is to base $2$. Also, we may assume without loss of generality that
$\omega(n)\leq n^{2}$, since otherwise just take $\omega'(n):=\min\{\omega(n),n^{2}\}$.
For $n\geq0$ denote the Taylor coefficients of $F$ by $c(n)$:\[
c(n):=\widehat{F}(n),\quad n\geq0\]
(remember that for $n<0$, $c(n)$ are defined by (\ref{eq:defcnneg})).

Next, define for $k\in\mathbb{N}$,
$$
f_{k}(z)   :=\sum_{n=-2^{k}}^{\infty}c(n)z^{n},   \quad z   \in\mathbb{D},\qquad
l_{k}(z)   :=\log|f_{k}(z)|,\qquad 
r_{k}   :=1-4^{-k}.
$$
A simple calculation shows that \begin{equation}
|f_{k}(z)-F(z)|\leq C\quad\forall|z|\geq1-2^{-k}.\label{eq:fkfconst}\end{equation}
Hence, from (\ref{eq:Fpoly}) we get 
\begin{align}\label{eq:lklinear}
F(z)=O\left(\frac{1}{(1-|z|)^{M}}\right) 
& \Rightarrow f_{k}=O\left(\frac{1}{(1-|z|)^{M}}\right)\textrm{ uniformly in }k \\
 & \Rightarrow l_{k}(z)\leq C_{2}k\quad\forall z,\:1-2^{-k}\leq\mathbf{|}z|\leq1-8^{-k}.
 \nonumber\end{align}
Define $A_{k}:=-\omega(k)/2$.

\begin{lem}
\label{lem:existsK}For every $\delta>0$ there exists a $K=K(\delta)$
such that \[
\mathbf{m}\{\theta:l_{K}(\theta r_{K})>A_{K}\}<\delta.\]

\end{lem}
\Proof 
Use Lemma \ref{lem:privdom} to find a compact $S\subset\partial\mathbb{D}$
of measure $>1-\delta/2$ satisfying the fact that $F$ is continuous on $\mathcal{P}(S)$.
Let $S'\subset S$ be a compact set of measure $>1-\delta$ such that\begin{equation}
\lim_{\eta\to0}\frac{\mathbf{m}(S\cap I(\theta,\eta))}{\eta}=1\textrm{ uniformly in }\theta\in S'.\label{eq:lebeguep}\end{equation}

Our purpose is to use Lemma \ref{lem:harmpriv} with $\varepsilon=1/4C_{1}$.
Therefore,  define  $\eta_{0}$ by the condition\[
\frac{\mathbf{m}(S\cap I(\theta,\eta))}{\eta}\geq1-\varepsilon^{2}\quad\forall\theta\in S',\forall\eta<\eta_{0}.\]
For any $k$ and any $z=\theta r_{k}$, $\theta\in S'$, let $\mathcal{A}=\mathcal{A}(k,\theta)$
be the component of the set $\mathcal{P}(S)\setminus(1-4^{-k}/\varepsilon)\mathbb{D}$
containing $z$. Examine the harmonic measure $\Omega$ of $\mathcal{A}$.
Lemma \ref{lem:harmpriv} gives, for all $k$ sufficiently large such
that $4^{-k}<\eta_{0}\varepsilon$, \begin{equation}
\Omega(z,\mathcal{A})(S\cap\partial\mathcal{A})\geq\frac{3}{4}.\label{eq:omegE34}\end{equation}

Now, if $k$ is sufficiently large so as to satisfy in addition that
$4^{-k}/\varepsilon<2^{-k}$ then we can use (\ref{eq:fkfconst}) and
the continuity of $F$ on $\mathcal{P}(S)$ to get \begin{equation}
l_{k}(z)\leq C\quad\forall k,\,\forall\theta\in S',\,\forall z\in\partial\mathcal{A}(k,\theta).\label{eq:lknotE}\end{equation}
 (\ref{eq:lknotE}) will be used on $\partial\mathcal{A}\setminus S$.
On $\partial\mathcal{A}\cap S$ we write \begin{equation}
|f_{k}(\theta)|\leq\sum_{n<-2^{k}}|c(n)|\leq C\sum_{n<-2^{k}}e^{-\omega(\log|n|)}\leq C\sum_{j=k}^{\infty}2^{j}e^{-\omega(j)}\leq Ce^{-(1-o(1))\omega(k)}\label{eq:fksmall}\end{equation}
(in the last inequality we used the fact  that $\omega(n)/n$ is increasing
to infinity, due to (\ref{eq:sumwfin})). Hence if $k$ is sufficiently
large we get\begin{equation}
l_{k}(\theta)\leq-\frac{3}{4}\omega(k)+C\quad\forall\theta\in S.\label{eq:lk34}\end{equation}
 Now $l$ is the logarithm of an analytic function and is therefore
subharmonic. This allows us  to use (\ref{eq:subkak}) and from (\ref{eq:omegE34}),
(\ref{eq:lknotE}) and (\ref{eq:lk34}) we get\[
l_{k}(z)\leq-\frac{9}{16}\omega(k)+C.\]
With $K$ sufficiently large  to satisfy all requirements
so far, as well as $\omega(K)/16\break >C$, the lemma is proved.
\Endproof\vskip4pt  

4.4.\quad  Continuing the proof of the theorem, let $\delta>0$ be some arbitrary
number, and let $z_{0}\in\mathbb{D}$ satisfy $|z_{0}|=1-\sqrt{\delta}$.
Let $P_{k}:=\Omega(z_{0},r_{k}\mathbb{D})$ (this is just the Poisson
kernel with appropriate parameters). Naturally, we assume $z_{0}\in r_{k}\mathbb{D}$,
so that everything below holds for $k>\log_{4}1/(1-|z_{0}|)$. Our
purpose is to show that the integrals \[
\int l_{k}(z)\, dP_{k}(z)\]
increase only in a precisely controlled manner. It turns out that
this is difficult to do directly, and we need to ``regularize''
before doing so. Define therefore\[
\mathcal{I}_{k}:=\int\left[l_{k}(z)\right]_{A_{k}}dP_{k}(z)\]
where $\left[f\right]_{M}:=\max(f,M)$. The proof will revolve around
a comparison of $\mathcal{I}_{k+1}$ and $\mathcal{I}_{k}$. Since
the calculation is long, we shall perform it in two stages, first
moving the circle of integration inward but keeping the $f_{k}$ and
only then changing $f_{k}$. 

\begin{lem}
\label{lem:LkprmLkp1}With the notation above{\rm ,}
\[
\int\left[l_{k+1}(z)\right]_{A_{k+1}}dP_{k}(z)
\leq\mathcal{I}_{k+1}+Ck^{2}2^{-k}.\]
\vskip8pt

\end{lem}
\Proof 
We wish to compare the harmonic measure on $r_{k+1}\mathbb{D}$ to
the one on the annulus $\mathcal{A}:=\{1-2^{-k-1}\leq|z|\leq r_{k+1}\}$.
The probability of a Brownian motion starting from $z$, $|z|=r_{k}$
to hit $\left\{ |z|=1-2^{-k-1}\right\} $ before $\left\{ |z|=r_{k+1}\right\} $
is $\leq C2^{-k}$ (this follows from the explicit formula for the
solution of Dirichlet's problem in an annulus (\ref{eq:harmann}))
, and therefore\begin{equation}
\left\Vert \Omega(z,\mathcal{A})-\Omega(z,r_{k+1}\mathbb{D})\right\Vert \leq C2^{-k}\quad\forall|z|=r_{k}\label{eq:QzQzp}\end{equation}
where the norm is the usual norm in the space of measures on $\mathbb{D}$.
We use this with the subharmonic function $\varphi(z):=\left[l_{k+1}(z)\right]_{A_{k+1}}$
and get\begin{align}\label{eq:lkAkz}
\varphi(z) & \stackrel{(\ref{eq:subkak})}{\leq}\int\varphi\, d\Omega(z,\mathcal{A})
\stackrel{(\ref{eq:QzQzp})}{\leq}\int\varphi\,
d\Omega(z,r_{k+1}\mathbb{D})+C2^{-k}\max_{\partial\mathcal{A}}|\varphi| \\
 & \stackrel{(\ref{eq:lklinear})}{\leq}\int\varphi\, d\Omega(z,r_{k+1}\mathbb{D})+C2^{-k}
(k+|A_{k+1}|) \nonumber \\
 & \,\leq\int\varphi\,
d\Omega(z,r_{k+1}\mathbb{D})+C2^{-k}k^{2}.\nonumber
\end{align} This we integrate
to get\begin{align*}
\int\left[l_{k+1}(z)\right]_{A_{k+1}}dP_{k}(z) & \stackrel{(\ref{eq:lkAkz})}{\leq}\int\int\left[l_{k+1}(x)\right]_{A_{k+1}}d\Omega(z,r_{k+1}\mathbb{D})(x)\, dP_{k}(z)+C2^{-k}k^{2}\\
 & =\int\left[l_{k+1}(z)\right]_{A_{k+1}}dP_{k+1}(z)+C2^{-k}k^{2}.\end{align*}
The last equality is the well known semigroup property of the Poisson
kernel. In probabilistic terms, this is integration over conditional
expectation.
\hfill\qed

\begin{lem}
\label{lem:epskepskp1}With the notation  above{\rm ,}
\begin{equation}
\mathcal{I}_{k}\leq\int\left[l_{k+1}(z)\right]_{A_{k}}dP_{k}(z)+Ce^{-c\omega(k)}.\label{eq:lemepskepskp}\end{equation}

\end{lem}
\Proof 
The difference between $f_{k}$ and $f_{k+1}$ can be estimated as
in (\ref{eq:fksmall})  (for $|z|=r_{k}$) by \[
|f_{k}(z)-f_{k+1}(z)|\leq\sum_{n=-2^{k+1}}^{-2^{k}-1}|c(n)|r_{k}^{n}\leq C\sum e^{-\omega(\log_{2}n)}\leq Ce^{-(1-o(1))\omega(k)}.\]
Therefore, if $|f_{k}(z)|\geq e^{-\omega(k)/2}$ then $|f_{k}|\leq|f_{k+1}|(1+Ce^{-\omega(k)(1/2-o(1))})$
which gives\[
l_{k}(z)\leq l_{k+1}(z)+Ce^{-c\omega(k)}\quad\forall|z|=r_{k},\, l_{k}(z)\geq A_{k}\]
or in other words\[
\left[l_{k}\right]_{A_{k}}\leq\left[l_{k+1}\right]_{A_{k}}+Ce^{-c\omega(k)}\]
which immediately gives (\ref{eq:lemepskepskp}) and proves the \pagebreak lemma.
\Endproof 
 
The gap between the $\left[l_{k+1}\right]_{A_{k+1}}$ of Lemma \ref{lem:LkprmLkp1}
and the $\left[l_{k+1}\right]_{A_{k}}$ of Lemma \ref{lem:epskepskp1}
is bridged by Lemma \ref{lem:main}. It will be convenient to reparametrize
as follows: define $B_{k}=C_{2}k$ so that by (\ref{eq:lklinear})
we would have \[
A_{k}\leq\left[l_{k}\right]_{A_{k}}(z)\leq B_{k}\quad\forall|z|=r_{k-1}.\]
Define $\varepsilon_{k}'$ by the relation\[
\int\left[l_{k}(z)\right]_{A_{k}}dP_{k-1}(z)=\varepsilon_{k}'B_{k}+(1-\varepsilon_{k}')A_{k}.\]
Then Lemma \ref{lem:main} gives
 \[
\int\left[l_{k+1}(z)\right]_{A_{k}}dP_{k}(z)
\leq\varepsilon_{k+1}'B_{k+1}+(1-\varepsilon_{k+1}')A_{k}.\]
Lemma \ref{lem:LkprmLkp1} and then Lemma \ref{lem:epskepskp1} now
give\begin{align}
\varepsilon_{k}'B_{k}+(1-\varepsilon_{k}')A_{k} & =\int\left[l_{k}(z)\right]_{A_{k}}dP_{k-1}(z)
\leq\mathcal{I}_{k}+Ck^{2}2^{-k}\label{eq:epspuenta} \\
 & \leq\varepsilon_{k+1}'B_{k+1}+(1-\varepsilon_{k+1}')A_{k}+Ck^{2}2^{-k}+Ce^{-c\omega(k)}
\nonumber\end{align}
so that
\begin{equation}
\varepsilon_{k}'\leq\varepsilon_{k+1}'\left(1+\frac{C_{2}}{B_{k}-A_{k}}\right)+\frac{Ck^{2}2^{-k}}{B_{k}-A_{k}}\leq\varepsilon_{k+1}'\left(1+\frac{C}{\omega(k)}\right)+Ck2^{-k}.\label{eq:epskepskp1}\end{equation}
The same holds for the more natural quantity $\varepsilon_{k}$ defined
by \begin{equation}
\int\left[l_{k}(z)\right]_{A_{k}}dP_{k}(z)=\varepsilon_{k}B_{k}+(1-\varepsilon_{k})A_{k}.\label{eq:defepsk}\end{equation}
Indeed, (\ref{eq:epspuenta}) shows that\[
\varepsilon_{k}'-Ck2^{-k}\leq\varepsilon_{k}\leq\varepsilon_{k+1}'\left(1+\frac{C}{\omega(k)}\right)+Ce^{-c\omega(k)}.\]

\demo{{\rm 4.5}} With this, the theorem is now easy. Lemma \ref{lem:existsK} combined
with the fact that $|P_{K}(z)|\leq C\delta^{-1/2}$ for all $z\in\partial(r_{K}\mathbb{D})$
shows that for some $K$ sufficiently large \[
\int\left[l_{K}(z)\right]_{A_{K}}dP_{K}(z)\leq CK\sqrt{\delta}-(1-C\sqrt{\delta})\frac{\omega(K)}{2}\]
so that $\varepsilon_{K}\leq C\sqrt{\delta}$. Applying (\ref{eq:epskepskp1})
repeatedly, we get that for all applicable $k$, $\varepsilon_{k}\leq C\sqrt{\delta}+Ck2^{-k}$
(here is where we use the condition $\sum\frac{1}{\omega(k)}<\infty$).
Let $k_{0}$ be the minimal applicable $k$, namely $\left\lceil \log_{2}1/\delta\right\rceil +1$.
In particular\begin{equation}
\varepsilon_{k_{0}}\leq C\sqrt{\delta}.\label{eq:epsk0}\end{equation}
 We can now estimate $l_{k_{0}}(z_{0})$ as in Lemma \ref{lem:LkprmLkp1}:
define $\mathcal{A}:=\{1-2^{-k_{0}}\leq|z|\leq r_{k_{0}}\}$ and get\begin{align*}
l_{k_{0}}(z_{0}) & \leq\int l_{k_{0}}\, d\Omega(z_{0},\mathcal{A})\stackrel{(\ref{eq:QzQzp})}{\leq}\int l_{k_{0}}\, dP_{k_{0}}+C2^{-k_{0}}\max_{\partial\mathcal{A}}|l_{k_{0}}|\\
 &
\!\stackrel{(*)}{\leq}\varepsilon_{k_{0}}B_{k_{0}}+ 
(1-\varepsilon_{k_{0}})A_{k_{0}}+Ck_{0}^{2}2^{-k_{0}}
\stackrel{(\ref{eq:epsk0})}{\leq}-c\omega(\log_{2}1/\delta) ,\end{align*}
where $(*)$ comes from the definition of $\varepsilon_{k}$, (\ref{eq:defepsk}) for the left term and
(\ref{eq:lklinear}) and $\omega(k)\leq k^{2}$ for the right term. Therefore $|f_{k_{0}}(z_{0})|\leq
e^{-c\omega(\log_{2}1/\delta)}$. Since this holds for all $z_{0}$ with $|z_{0}|=1-\sqrt{\delta}$
we may estimate $c(n)$ using Laurent's formula (we need to assume
$n\geq-2^{k_{0}}$ so assume $n\geq-1/\delta$) and get \[
|c(n)|\leq(1-\sqrt{\delta})^{-1-n}e^{-c\omega(\log_{2}1/\delta)}.\]
Since this holds for all $\delta>0$, we get that $c(n)=0$ for all
$n$, and the theorem is proved.\hfill\qed

\end{document}